\documentclass[11pt]{article}
\usepackage{amsmath,amssymb,amscd,enumerate,eucal}
\numberwithin{equation}{section}
\newcommand{\PP}{\mathcal{P}}

\newcommand{\R}{\mathbb{R}}

\newcommand{\Z}{\mathbb{Z}}

\newcommand{\Ad}{\text{Ad}}

\newcommand{\HH}{\mathbb{H}}
\newcommand{\YMe}{\mathcal{YM}_{\epsilon}}

\newcommand{\q}{\mathsf{q}}
\newcommand{\A}{\mathsf{a}}

\newcommand{\Je}{\mathcal{J}_{\epsilon}}
\newcommand{\tr}{\mathsf{Tr}}
\newcommand{\so}{\mathfrak{so}}
\newcommand{\cat}{\text{cat}}
\newcommand{\Asf}{\mathsf{A}}
\newcommand{\hsf}{\mathsf{h}}

\begin{document}
\setlength{\baselineskip}{17pt}

\title{Small coupling limit and multiple solutions to the Dirichlet
Problem for Yang-Mills connections in $4$ dimensions - Part II}
\author{Takeshi Isobe \thanks {%
Tokyo Institute of Technology; email:
isobe.t.ab@m.titech.ac.jp} $\,$ and Antonella Marini \thanks {%
University of L'Aquila / Yeshiva University; email:
marini@dm.univaq.it}}


\date{}
\maketitle

\begin{abstract}
{\footnotesize \baselineskip 4mm} In this paper we complete the
proof of the existence of multiple solutions (and, in particular,
non minimal ones), to the $\epsilon$-Dirichlet problem obtained as a
variational problem for the $SU(2)_\epsilon$-Yang Mills functional.
This is equivalent to proving the existence of multiple solutions to
the Dirichlet problem for the $SU(2)$-Yang Mills functional with
small boundary data. In the first paper of this series this
non-compact variational problem is transformed into the finite
dimensional problem of finding the critical points of the function
$\Je(\q)$, which is essentially the Yang Mills functional evaluated
on the approximate solutions, constructed via a gluing technique. In
the present paper, we establish a Morse theory for $\Je(\q)$, by
means of Ljusternik-Schnirelmann theory, thus complete the proofs of
Theorems 1-3 in \cite{IM1}.
\end{abstract}



\section{Introduction}

Let $\mathcal{A}_{+1}(A_0)$ be space of smooth connections with relative 2nd Chern number equal to $1$, calculated with respect to a fixed absolute minimizer, say $\underline{A}_{\epsilon}$, with prescribed smooth boundary value
$A_0,$ of the
$SU(2)_\epsilon$-Yang Mills functional
\begin{equation}
\label{YMe} \YMe(A)=\int_{B^4}|{F_A}^\epsilon|^2\,dx\,,
\end{equation}
where
${F_A}^\epsilon=dA+\frac{1}{2}[A,A]_{\epsilon}:=dA+\frac{\epsilon}{2}[A,A].$ The space $\mathcal{A}_{+1}(A_0)$ is well-defined and non-empty, since an absolute minimum is known to exist always (cf. \cite{Marini}), and by Taubes' gluing construction (cf. \cite{Taubes}).

The Euler Lagrange equations for \eqref{YMe} with prescribed smooth boundary value $A_0$ yield the
$\epsilon$-Dirichlet problem:
$$
\bigl(\mathcal{D}_{\epsilon}\bigr)\quad\qquad\left\{\begin{array}{ll}
{d_A^{\ast}}^\epsilon {F_A}^\epsilon=0\quad&\mbox{in } B^4\\
\iota^{\ast}A\sim A_0\quad&\mbox{at } \partial B^4\,,
\end{array}\right.
$$
where, $\iota:\partial B^4\to\overline{B}^4$ is the inclusion, the
symbol $\sim$  stands for gauge equivalence,
${d_A^{\ast}}^\epsilon:=\ast d\ast+\epsilon\ast[A,\ast\cdot]$,  and
$\ast$ denotes the Hodge star operator with respect to the flat
metric on $\R^4$.

We recall that, for $\epsilon > 0$, solutions to
$\bigl(\mathcal{D}_{\epsilon}\bigr)$ are in one-to-one
correspondence with solutions to
$\bigl(\mathcal{D}(\epsilon)\bigr)$, that is the Dirichlet problem
for the standard $SU(2)$-Yang Mills functional with boundary value
$\epsilon A_0$ (cf. $\S$ 2.2 in \cite{IM1}).

In \cite{IM} it is shown that a solution to
$(\mathcal{D}_{\epsilon})$, i.e. a \emph{large} Yang Mills field
with boundary value $A_0$, always exists in $\mathcal{A}_{+1}(A_0)$,
or $\mathcal{A}_{-1}(A_0)$ if $A_0$ is non-flat. With no loss of
generality, we focus on $\mathcal{A}_{+1}(A_0)$ (since all the
results obtained can be transformed into results on
$\mathcal{A}_{-1}(A_0)$, by simply reversing the orientation of
$B^4).$

In the first paper of this series, the problem of finding solutions to the $\epsilon$-Dirichlet problem
$(\mathcal{D}_{\epsilon})$ in
$\mathcal{A}_{+1}(A_0)$ is reduced, for small values of the
parameter $\epsilon>0$,  to the finite dimensional problem of
finding critical points of
\begin{equation}
\label{curl}
\Je(\q)=\epsilon^2\int_{B^4}|{F_{A(\q)+a(\q)}}^{\epsilon}|^2\,dx\,, \quad\mbox{ for }\q\in\PP(d_0,\lambda_0;D_1,D_2,\epsilon)\,,
\quad \mbox{ with } 0<D_1<D_2\,,
\end{equation}
where
\begin{equation}
\label{subpar}
\PP(d_0,\lambda_0;D_1,D_2;\epsilon):=\{\q:=(p,[g],\lambda)\in
\PP(d_0,\lambda_0):D_1\epsilon<\lambda^2<D_2\epsilon\}\,,
\end{equation}
with
\begin{equation}
\label{par}
\PP(d_0,\lambda_0):= B^4_{1-d_0}\times SO(3)\times(0,\lambda_0)\,, \quad\mbox{ for }0<2\lambda_0<d_0\,,
\end{equation}
is the parameter space used in the gluing procedure to construct the approximate solutions $A(\q)$.

We recall that $A(\q)$ is the approximate solution obtained by
gluing the $1$-instanton to  $\underline{A}_{\epsilon}$, and $a(\q)$
is the unique solution to the \emph{auxiliary equation} defined in
$\S3.6$ of \cite{IM1} (cf., in particular, Lemma 3.9 in \cite{IM1}),
and accounts for the interaction with the boundary (cf. \cite{IM1}
for the precise definitions). In Proposition 3.2 of \cite{IM1} it is
proved that the critical points of $\Je(\q)$ are in one-to-one
correspondence with the solutions to ${\mathcal D}_\epsilon.$ In the
present paper, we do the Morse theory for the finite dimensional
functional $\Je(\q)$ on $\PP(d_0,\lambda_0;D_1,D_2;\epsilon)$ and
complete the proof of the existence theorems of multiple solutions
to ${\mathcal D}_\epsilon$ (Theorems \ref{1}-\ref{3} in $\S5-\S7$).

\medskip\noindent
More in detail, in $\S2$ we study the asymptotic properties, as $p$
approaches $\partial B^4$, of the function $F(p)$ and of the
$3\times 3$ matrix $M(A_0,p)$, which play a crucial role in our
construction; in $\S3$ we study the Morse theoretical properties of
the function $R\to\tr(RM),$ with $R\in SO(3),$ which translate into
properties of the space of solutions to $({\mathcal D}_\epsilon)$,
by taking $M=M(A_0,p)$; in $\S4$ we prove two lemmas on the
asymptotic behavior, as $\epsilon\to 0$, of the functional
$\Je(\q)$, needed for the proofs of the main theorems, Theorems
\ref{1}-\ref{3}, proved in $\S5-\S7$; in $\S8$ we give a method to
construct boundary data that yield any given matrix $M= M(A_0, p)$
(cf. next section), thus all the different cases enlisted in Theorem
\ref{3} can be obtained. Moreover, we show that the non-degeneracy
assumption $\mu_1\ge\mu_2\ge \mu_3> 0$ on the eigenvalues $\mu_i$
($i =1,2,3$) of $M(A_0, p)^t M(A_0, p)$ can always be achieved via a
small perturbation of the boundary data.

\medskip\noindent
We assume throughout this paper that $\lambda_0,$ $d_0,$ $D_1,$
$D_2,$ $\epsilon$, $\q$ in \eqref{curl}--\eqref{par} be fixed once
and for all,  with $0<D_1<D_2,$ and $0<2\lambda_0<d_0.$

\section{Behavior of $F(p)$ and $M(A_0,p)$ near the boundary}

The existence of multiple, and in particular non-minimal solutions
to the $\epsilon$-Dirichlet problem $(\mathcal{D}_{\epsilon})$ for
small positive values of the parameter $\epsilon$, depends
critically on the properties of the function $F(p)$ and of the
$3\times 3$ matrix $M(A_0,p)$, as expressed in the statements of
Theorems \ref{1}-\ref{3} (cf. $\S$5-$\S$7).

\noindent We recall that
\begin{equation}
\label{F}
F(p)=\int_{B^4}|(dh_p)^-|^2\,dx\;,\; p\in B^4\;,
\end{equation}
\noindent
and that, for a given boundary value $A_0$,
\begin{equation}
\label{M}
M(A_0,p):=\bigl(m_{ij}(A_0,p)\bigr)\;,
\end{equation}
with
$$m_{ij}(A_0,p):=\int_{B^4}\bigl((d\underline{A}_{0,j})^-,(dh_{p,i})^-\bigr) \,dx\;,\qquad (1\le i,j\le 3)\,,$$
where, for a given $2$-form $\omega$, we denote its anti-self dual
component by $\omega^-$ (i.e., $\omega^-:=(\omega-\ast\omega)/2$),
and $\underline{A}_0$ is a solution to the linear Dirichlet problem
$$\bigl(\mathcal{D}_0\bigr)\qquad\quad\left\{\begin{array}{ll}
d^{\ast}dA=0\quad&\mbox{in } B^4\\
\iota^{\ast}A\sim A_0\quad&\mbox{on } \partial B^4\;,
\end{array}\right.$$
(Note that, by Hodge theory, $d \underline{A}_0$ is uniquely
determined by the boundary value $A_0$, thus the definition above is
well-posed).

\noindent In this section we study the asymptotic behavior of $F(p)$ and $M(A_0,p)$, as $p$
approaches the boundary $\partial B^4$ of the four dimensional disk.

\newtheorem{lemma}{Lemma}[section]
\begin{lemma}
\label{L4.3} The function $F(p)=\int_{B^4}|(dh_p)^-|^2\,dx$, for $p\in B^4$ satisfies:
\begin{enumerate}[(1)]
\item $F(p)>0$ for all $p\in B^4$;
\item there exists a constant $C_1>0$ (independent of $p$) such that $F(p)=C_1d(p)^{-4}+o(d(p)^{-4})$ as $p\in B^4$ approaches $\partial B^4$,
where
$d(p):=1-|p|$ is the distance from $p$ to $\partial B^4$;
\item there exists a constant $C_2>0$ (independent of $p$) such that $F'(p)=C_2d(p)^{-5}\frac{p}{|p|}+o(d(p)^{-5})$
as $p\in B^4$ approaches $\partial B^4$.
\end{enumerate}
\end{lemma}
\textit{Proof:}

\noindent\emph{Proof of (1).}  Recall that $h_p$ is defined as the
solution of the Dirichlet problem $\Delta h_p=0$ in $B^4$ with boundary data
$h_p(x)=\text{Im}\frac{(\overline{x}-\overline{p})dx}{|x-p|^4}$ at
$\partial B^4$ (where all the components, not only the tangential ones, are assigned at the boundary).
For $1\le i\le 4$, let
$\alpha_{p,i}=\alpha_{p,i}(x)$ be the solution of the Dirichlet
problem $\Delta \alpha_{p,i}=0$ in $B^4$ and
$\alpha_{p,i}(x)=\frac{x_i-p_i}{|x-p|^4}$ at $\partial B^4$. To
prove (1), we rely on the Poisson integral representation of these
functions.

\noindent By the Poisson's formula (see~\cite{Taylor}), one has
\begin{equation}
\label{4.13}
\alpha_{p,i}(x)=\frac{1-|x|^2}{2\pi^2}\int_{S^3}\frac{1}{|x-y|^4}\frac{y_i-p_i}{|y-p|^4}\,dy
\end{equation}
and, writing $h_p=h_{p,1}i+h_{p,2}j+h_{p,3}k$, these components
satisfy
\begin{align}
\label{4.14}
&h_{p,1}=-\alpha_{p,2}dx^1+\alpha_{p,1}dx^2+\alpha_{p,4}dx^3-\alpha_{p,3}dx^4, \notag\\
&h_{p,2}=-\alpha_{p,3}dx^1-\alpha_{p,4}dx^2+\alpha_{p,1}dx^3+\alpha_{p,2}dx^4, \notag\\
&h_{p,3}=-\alpha_{p,4}dx^1+\alpha_{p,3}dx^2-\alpha_{p,2}dx^3+\alpha_{p,1}dx^4.
\end{align}
By direct computation,
\begin{align}
\label{4.17}
(dh_{p,1})^-&=\frac{1}{2}\Big(\frac{\partial\alpha_{p,1}}{\partial x^1}+\frac{\partial\alpha_{p,2}}{\partial x^2}+
\frac{\partial\alpha_{p,3}}{\partial x^3}+\frac{\partial\alpha_{p,4}}{\partial x^4}\Big)\omega_1^-\notag\\
&\quad+\frac{1}{2}\Big(-\frac{\partial\alpha_{p,1}}{\partial x^4}+\frac{\partial\alpha_{p,2}}{\partial x^3}-
\frac{\partial\alpha_{p,3}}{\partial x^2}+\frac{\partial\alpha_{p,4}}{\partial x^1}\Big)\omega_2^-\notag\\
&\quad+\frac{1}{2}\Big(\frac{\partial\alpha_{p,1}}{\partial
x^3}+\frac{\partial\alpha_{p,2}}{\partial
x^4}-\frac{\partial\alpha_{p,3}}{\partial
x^1}-\frac{\partial\alpha_{p,4}}{\partial x^2}\Big)\omega_3^-,
\end{align}
\begin{align}
\label{4.18}
(dh_{p,2})^-&=\frac{1}{2}\Big(\frac{\partial\alpha_{p,1}}{\partial x^4}-\frac{\partial\alpha_{p,2}}{\partial x^3}+
\frac{\partial\alpha_{p,3}}{\partial x^2}-\frac{\partial\alpha_{p,4}}{\partial x^1}\Big)\omega_1^-\notag\\
&\quad+\frac{1}{2}\Big(\frac{\partial\alpha_{p,1}}{\partial x^1}+\frac{\partial\alpha_{p,2}}{\partial x^2}+\frac{\partial\alpha_{p,3}}{\partial x^3}+
\frac{\partial\alpha_{p,4}}{\partial x^4}\Big)\omega_2^-\notag\\
&\quad+\frac{1}{2}\Big(-\frac{\partial\alpha_{p,1}}{\partial
x^2}+\frac{\partial\alpha_{p,2}}{\partial
x^1}+\frac{\partial\alpha_{p,3}}{\partial x^4}-
\frac{\partial\alpha_{p,4}}{\partial x^3}\Big)\omega_3^-,
\end{align}
and
\begin{align}
\label{4.19}
(dh_{p,3})^-&=\frac{1}{2}\Big(-\frac{\partial\alpha_{p,1}}{\partial
x^3}- \frac{\partial\alpha_{p,2}}{\partial
x^4}+\frac{\partial\alpha_{p,3}}{\partial x^1}+
\frac{\partial\alpha_{p,4}}{\partial x^2}\Big)\omega_1^-\notag\\
&\quad+\frac{1}{2}\Big(\frac{\partial\alpha_{p,1}}{\partial x^2}-\frac{\partial\alpha_{p,2}}{\partial x^1}-
\frac{\partial\alpha_{p,3}}{\partial x^4}+
\frac{\partial\alpha_{p,4}}{\partial x^3}\Big)\omega_2^-\notag\\
&\quad+\frac{1}{2}\Big(\frac{\partial\alpha_{p,1}}{\partial
x^1}+\frac{\partial\alpha_{p,2}}{\partial
x^2}+\frac{\partial\alpha_{p,3}}{\partial
x^3}+\frac{\partial\alpha_{p,4}}{\partial x^4}\Big)\omega_3^-,
\end{align}
where $\omega_1^-:=dx^1\wedge dx^2-dx^3\wedge dx^4$,
$\omega_2^-:=dx^1\wedge dx^3+dx^2\wedge dx^4$ and
$\omega_3^-:=dx^1\wedge dx^4-dx^2\wedge dx^3$ compose the standard basis
for anti-self dual forms on $\R^4$.

\noindent
By \eqref{4.13}, we calculate
\begin{equation}
\label{4.20}
\frac{\partial\alpha_{p,1}}{\partial x^1}(0)+\frac{\partial\alpha_{p,2}}{\partial x^2}(0)+
\frac{\partial\alpha_{p,3}}{\partial x^3}(0)+\frac{\partial\alpha_{p,4}}{\partial x^4}(0)=
\frac{2}{\pi^2}\int_{S^3}\frac{1-y\cdot p}{|y-p|^6}\,dy>0
\end{equation}
for all $p\in B^4$, where ``$\cdot$" represents the inner product in
$\R^4$.

From \eqref{4.17}--\eqref{4.19}, it follows that
$|(dh_{p,i})^-(0)|>0$ for $1\le i\le 3$ and for all $p\in B^4$, thus
$F(p)=\int_{B^4}|(dh_p)^-|^2\,dx>0$ for all $p\in B^4$.

\medskip

\noindent\emph{Proof of (2).} To prove (2), we write the functions
$\alpha_{p,i}(x)$ explicitly in terms of the Green function for the
Laplacian on $B^4$. More precisely, let
$\Gamma(x,y)=\Gamma(|x-y|)=-\frac{1}{8\omega_4}|x-y|^{-2},$  where $\omega_4$ is the volume of $B^4$, be the
fundamental solution for the Laplacian in $\R^4$, and $G(x,y)$ the
Green's function on $B^4$, with Dirichlet boundary
data. Denote by
$H(x,y)$ the regular part of $G(x,y)$, i.e.,
$H(x,y)=\Gamma(x,y)-G(x,y)$. The functions $\alpha_{p,i}$ are then
given by
\begin{equation}
\notag
\alpha_{p,i}(x)=-4\omega_4\frac{\partial H}{\partial p_i}(x,p)\,,
\end{equation}
and $H(x,p)$ is explicitly given by
$H(x,p)=\Gamma(|p||x-p_{\ast}|)$ (c.f.~\cite{GT}), where
$p_{\ast}=p/|p|^2$ for $p\ne 0$, and $p_{\ast}=\infty$ for $p=0$.
We thus have
\begin{equation}
\notag
\alpha_{p,i}(x)=-\frac{p_i}{|p|^4|x-p_{\ast}|^2}+\frac{x_i-p_{\ast,i}}{|p|^4|x-p_{\ast}|^4}-
\frac{2p_ip\cdot(x-p_{\ast})}{|p|^6|x-p_{\ast}|^4}\,,
\end{equation}
and finally
\begin{align}
\label{4.23}
\frac{\partial\alpha_{p,i}}{\partial x_j}=&\,\frac{2p_i(x_j-p_{\ast,j})}{|p|^4|x-p_{\ast}|^4}+
\frac{\delta_{ij}}{|p|^4|x-p_{\ast}|^4}-\frac{4(x_i-p_{\ast,i})(x_j-p_{\ast,j})}{|p|^4|x-p_{\ast}|^6}\notag\\
&-\frac{2p_ip_j}{|p|^6|x-p_{\ast}|^4}+\frac{8p_i(x_j-p_{\ast,j})p\cdot(x-p_{\ast})}{|p|^6|x-p_{\ast}|^6}.
\end{align}
Using \eqref{4.23}, one calculates explicitly \eqref{4.17}--\eqref{4.19}.
In particular, the coefficient of $\omega_1^-$ in \eqref{4.17}, namely
$((dh_{p,1})^-,\omega_1^-)$, is
\begin{equation}
\notag
((dh_{p,1})^-,\omega_1^-)=\frac{2p\cdot(x-p_{\ast})}{|p|^4|x-p_{\ast}|^4}-\frac{2}{|p|^4|x-p_{\ast}|^4}+
\frac{8(p\cdot(x-p_{\ast}))^2}{|p|^6|x-p_{\ast}|^6}\,.
\end{equation}
Similarly, the coefficients of $\omega_2^-$, $\omega_3^-$ in \eqref{4.17}, namely $((dh_{p,1})^-,\omega_2^-)$, $((dh_{p,1})^-,\omega_3^-)$, are
\begin{align}
\notag
((dh_{p,1})^-,\omega_2^-)&=\,\biggl(\frac{2}{|p|^4|x-p_{\ast}|^4}+
\frac{8p\cdot(x-p_{\ast})}{|p|^6|x-p_{\ast}|^6}\biggr)(p_2(x_3-p_{\ast,3})-p_3(x_2-p_{\ast,2}))\notag\\
&\quad+\biggl(\frac{2}{|p|^4|x-p_{\ast}|^4}+\frac{8p\cdot(x-p_{\ast})}{|p|^6|x-p_{\ast}|^6}\biggr)
(p_4(x_1-p_{\ast,1})-p_1(x_4-p_{\ast,4}))\,,
\end{align}
\begin{align}
\notag
((dh_{p,1})^-,\omega_3^-)&=\,\biggl(\frac{2}{|p|^4|x-p_{\ast}|^4}+\frac{8p\cdot(x-p_{\ast})}{|p|^6|x-p_{\ast}|^6}\biggr)
(p_2(x_4-p_{\ast,4})-p_4(x_2-p_{\ast,2}))\notag\\
&\quad+\biggl(\frac{2}{|p|^4|x-p_{\ast}|^4}+\frac{8p\cdot(x-p_{\ast})}{|p|^6|x-p_{\ast}|^6}\biggr)
(p_1(x_3-p_{\ast,3})-p_3(x_1-p_{\ast,1}))\,.
\end{align}
In order to study the behavior of the integral
$\int_{B^4}|(dh_p)^-|^2\,dx$ as $p\to\partial B^4$, it is sufficient
to take $p=(0,0,0,-1+d)$ and let $d\to 0$. (The general case follows from this since approaching from
tangential directions would contributes only lower order terms).  By taking
only the leading terms in the integrals below, one obtains the following asymptotic behaviors
as $d\to 0$ (notice that $p_{\ast}=(0,0,0,-1-d)+O(d^2)$ as:
$d\to 0$):
\begin{align}
\label{4.27}
&\int_{B^4}((dh_{p,1})^-,\omega_1^-)^2
\simeq\,\int_{B^4}\biggl(\frac{1}{|p|^4|x-p_{\ast}|^4}-\frac{4(p\cdot(x-p_{\ast}))^2}{|p|^6|x-p_{\ast}|^6}\biggr)^2\,dx
\notag\\
\simeq&\,\int_{x_4\ge -1}\biggl(\frac{1}{(x_1^2+x_2^2+x_3^2+(x_4+1+d)^2)^2}-
\frac{4(x_4+1+d)^2}{(x_1^2+x_2^2+x_3^2+(x_4+1+d)^2)^3}\biggr)^2\,dx\notag\\
\simeq&\,\frac{1}{d^4}\int_{x_4\ge 0}\biggl(\frac{1}{(x_1^2+x_2^2+x_3^2+(x_4+1)^2)^2}-
\frac{4(x_1+1)^2}{(x_1^2+x_2^2+x_3^2+(x_4+1)^2)^3}\biggr)^2\,dx
\simeq\,\frac{1}{d^4}\,.
\end{align}
Similarly,
\begin{align}
\label{4.28}
&\int_{B^4}((dh_{p,1})^-,\omega_2^-)^2 \notag\\
\simeq&\,\int_{B^4}\frac{(p\cdot(x-p_{\ast}))^2}{|p|^{12}|x-p_{\ast}|^{12}}(p_2(x_3-p_{\ast,3})-
p_3(x_2-p_{\ast,2})+p_4(x_1-p_{\ast,1})-p_1(x_4-p_{\ast,4}))^2\,dx\notag\\
\simeq&\,\int_{x_4\ge -1}\frac{(x_4+1+d)^2x_1^2}{(x_1^2+x_2^2+x_3^2+(x_4+1+d)^2)^6}\,dx
\simeq\,\frac{1}{d^4}\,,
\end{align}
\begin{align}
\label{4.29}
&\int_{B^4}((dh_{p,1})^-,\omega_3^-)^2\,dx\notag\\
\simeq&\,\int_{B^4}\frac{(p\cdot(x-p_{\ast}))^2}{|p|^{12}|x-p_{\ast}|^{12}}(p_2(x_4-p_{\ast,4})-p_4(x_2-p_{\ast,2})+
p_1(x_3-p_{\ast,3})-p_3(x_1-p_{\ast,1}))^2\,dx\notag\\
\simeq&\,\int_{x_4\ge-1}\frac{(x_4+1+d)^2x_2^2}{(x_1^2+x_2^2+x_3^2+(x_4+1+d)^2)^6}\,dx
\simeq\,\frac{1}{d^4}\,.
\end{align}
Combining \eqref{4.27} -- \eqref{4.29} yields
\begin{equation}
\int_{B^4}((dh_{p,j})^-)^2\,dx\simeq\frac{1}{d^4},\quad(d\to 0),\qquad\mbox{ for }\, j=1,2,3\,,
\end{equation}
thus
$\int_{B^4}|(dh_p)^-|^2\,dx\simeq\frac{1}{d^4}$ as $d\to 0$.

\medskip

\noindent \emph{Proof of (3).}

The proof above just showed that there exists a positive constant
$C$, independent of $p$, such that
$\int_{B^4}|(dh_p)^-|^2\,dx=Cd(p_{\ast},\partial
B^4)^{-4}+o(d(p_{\ast},\partial B^4)^{-4})$ as $d\to 0$. Formally,
the assertion in (3) follows from differentiating this equation with
respect to $p$. The detailed calculation, although more involved
than the one performed to prove (2), is conceptually the same.
So, we omit the details. \hfill$\Box$

\medskip

Next, we estimate the limit behavior of $M(A_0,p)$ as $p$ approaches
$\partial B^4$:

\begin{lemma}
\label{L4.4} Let $M(A_0,p):=(m_{ij}(A_0,p))$, for $p\in B^4$, be the
matrix defined in the introduction,  i.e.,
$m_{ij}(A_0,p)=\int_{B^4}\bigl((d\underline{A}_{0,j})^-,(dh_{p,i})^-\bigr)\,dx$.
Then $m_{ij}(A_0,p)$ and its derivatives $\frac{\partial
m_{ij}(A_0,p)}{\partial p_k}$ ($1\le k\le 4$) remain bounded as $p$
approaches the boundary.
\end{lemma}
\textit{Proof:} The harmonicity of $d\underline{A}_{0,j}$ and Stokes
theorem imply
\begin{equation}
\label{4.33} m_{ij}(A_0,p)=-\int_{\partial
B^4}\iota^{\ast}(d\underline{A}_{0,j})^-\wedge\iota^{\ast}(h_{p,i})\,
\end{equation}
where $\iota^{\ast}$ denotes the restriction to the boundary (the
pull-back via the inclusion). For simplicity, we write
$(d\underline{A}_{0,j})^-=\alpha_1\omega_1^-+\alpha_2\omega_2^-+\alpha_3\omega_3^-$
and, with no loss of generality, consider the case $i=1$. Writing
$\iota^{\ast}h_{p,1}=\frac{1}{|x-p|^4}\bigl(-(x_2-p_2)dx_1+(x_1-p_1)dx_2-(x_3-p_3)dx_4+(x_4-p_4)dx_3\bigr)$
at $\partial B^4$, we have explicitly
\begin{align}
\label{4.34}
m_{ij}(A_0,p)=&\,\int_{\partial B^4}\frac{1}{|x-p|^4}\bigl((\alpha_2(x)(x_1-p_1)-\alpha_3(x)(x_2-p_2)-\alpha_1(x)(x_4-p_4))dx_1\wedge dx_2\wedge dx_3\notag\\
&+(\alpha_1(x)(x_3-p_3)+\alpha_2(x)(x_2-p_2)+\alpha_3(x)(x_1-p_1))dx_1\wedge dx_2\wedge dx_4\notag\\
&+(-\alpha_1(x)(x_2-p_2)+\alpha_2(x)(x_3-p_3)+\alpha_3(x)(x_4-p_4))dx_1\wedge dx_3\wedge dx_4\notag\\
&+(\alpha_1(x)(x_1-p_1)-\alpha_3(x)(x_3-p_3)+\alpha_2(x)(x_4-p_4))dx_2\wedge dx_3\wedge dx_4\bigr)\notag\\
=&\,\int_{\partial B^4\cap|x-e|<\delta}\cdots+\int_{\partial B^4\cap|x-e|\ge\delta}\cdots,
\end{align}
where $e=p/|p|$.

As for the proof of Lemma \ref{L4.3}, it is sufficient to consider the case
$p=(0,0,0,-1+d)$, and let $d\to 0$. Under this assumption,
$e=(0,0,0,-1)$ and the second integral is bounded as $d\to 0$, for
any fixed $\delta>0$. For small positive $\delta$, the first
integral can be rewritten as an integral over
$B^3_{\sqrt{\delta^2-d^2}}$, via the coordinate transformation
$B^3_{\sqrt{\delta^2-d^2}}\ni x':=(x_1,x_2,x_3)\mapsto
(x',-\sqrt{1-|x'|^2})\in\partial B^4\cap\{|x-e|<\delta\}$. Since
$dx_4=O(|x|)$, it is easy to see that all the
integrals over $\partial B^4\cap\{|x-e|<\delta\}$ in \eqref{4.34} remain
bounded as $d\to 0$, with the possible exception of
\begin{equation}
\label{4.35}
\int_{\partial
B^4\cap|x-e|^2<\delta}\frac{1}{|x-p|^4}\bigl((\alpha_2(x)(x_1-p_1)-\alpha_3(x)(x_2-p_2)-\alpha_1(x)(x_4-p_4))dx_1\wedge
dx_2\wedge dx_3\bigr)
\end{equation}
To check that the latter also remains bounded as $d\to 0$, we write
(using Taylor's expansion)
\begin{align}
\label{4.36}
&\int_{\partial B^4\cap|x-e|^2<\delta}\frac{1}{|x-p|^4}\big((\alpha_2(x)(x_1-p_1)-\alpha_3(x)(x_2-p_2)-
\alpha_1(x)(x_4-p_4))dx_1\wedge dx_2\wedge dx_3\big)\notag\\
&=\int_{|x'|^2<\delta^2-d^2}\frac{1}{(|x'|^2+d^2)}(\alpha_2(x',-\sqrt{1-|x'|^2})x_1-\alpha_3(x',-\sqrt{1-|x'|^2})x_2\notag\\
&\quad-\alpha_1(x',-\sqrt{1-|x'|^2})(-\sqrt{1-|x'|^2}+1-d)\,dx'\notag\\
&=\int_{|x'|^2<\delta^2-d^2}\frac{1}{(|x'|^2+d^2)^2}(\alpha_2(e)x_1-\alpha_3(e)x_2+d\alpha_1(e))\,dx'\notag\\
&\quad+
O\biggl(\int_{|x'|^2<\delta^2-d^2}\frac{|x'|^2}{(|x'|^2+d^2)^2}\,dx'\biggr)+
O\biggl(\int_{|x'|^2<\delta^2-d^2}\frac{d|x'|}{(|x'|^2+d^2)^2}\,dx'\biggl)\notag\\
&=d\alpha_1(e)\int_{|x'|^2<\delta^2-d^2}\frac{1}{(|x'|^2+d^2)^2}\,dx'+
O\biggl(\int_{|x'|^2<\delta^2-d^2}\frac{|x'|^2}{(|x'|^2+d^2)^2}\,dx'\biggr)\notag\\
&\quad
+O\biggl(\int_{|x'|^2<\delta^2-d^2}\frac{d|x'|}{(|x'|^2+d^2)^2}\,dx'\biggl)=O(1)\qquad
(d\to 0).
\end{align}
Thus, $m_{ij}(A_0,p)$ stays bounded as $p$ approaches $\partial
B^4$.

Next we differentiate \eqref{4.33} with respect to  $p_k$, and obtain
\begin{equation}
\label{4.37}
\frac{\partial m_{ij}(A_0,p)}{\partial p_k}= -\int_{\partial
B^4}\iota^{\ast}({F_{\underline{A}_0,j}}^0)^-\wedge\iota^{\ast}\Big(\frac{\partial(
h_{p,i})^-}{\partial p_k}\Big).
\end{equation}
With no loss of generality, we take $k=1$. Again, the only troublesome contribution (as $d\to 0$) in \eqref{4.37}
could come from the
integral over $\partial B^4\cap\{|x-e|<\delta\}$ containing
$dx_1\wedge dx_2\wedge dx_3$. The term of
$-\frac{\partial}{\partial p_1}(F_{\underline{A}_0}^-,(dh_{p,1})^-)$
containing $dx_1\wedge dx_2\wedge dx_3$ is
\begin{align*}
&\frac{\partial }{\partial p_1}(|x-p|^{-4}(\alpha_2(x)(x_1-p_1)-\alpha_3(x)(x_2-p_2)-\alpha_1(x)(x_4-p_4))dx_1\wedge dx_2\wedge dx_3\notag\\
&=4|x-p|^{-6}(\alpha_2(x)(x_1-p_1)^2-\alpha_3(x)(x_1-p_1)(x_2-p_2)-\alpha_1(x)(x_1-p_1)(x_4-p_4))dx_1\wedge dx_2\wedge dx_3\notag\\
&\quad-|x-p|^{-4}\alpha_2(x)dx_1\wedge dx_2\wedge dx_3\,,
\end{align*}
and the corresponding integral over $\partial B^4\cap\{|x-e|<\delta\}$
is
\begin{align}
\label{4.38}
 &\int_{|x'|^2<\delta^2-d^2}\frac{4}{(|x'|^2+d^2)^3}(\alpha_2(x',-\sqrt{1-|x'|^2})x_1^2-
 \alpha_3(x',-\sqrt{1-|x'|^2})x_1x_2\notag\\
 &\quad -\alpha_1(x',-\sqrt{1-|x'|^2})x_1(-\sqrt{1-|x'|^2}+1-d)\,dx_1\wedge dx_2\wedge dx_3\notag\\
 &\quad-\int_{|x'|^2<\delta^2-d^2}\frac{\alpha_2(x',-\sqrt{1-|x'|^2})}{(|x'|^2+d^2)^2}\,dx'\notag\\
&=\int_{|x'|^2<\delta^2-d^2}\frac{4\alpha_2(e)x_1^2}{(|x'|^2+d^2)^3}\,dx'-
\int_{|x'|^2<\delta^2-d^2}\frac{4\alpha_3(e)x_1x_2}{(|x'|^2+d^2)^3}\,dx'+
\int_{|x'|^2<\delta^2-d^2}\frac{d\alpha_1(e)x_1}{(|x'|^2+d^2)^3}\,dx' \notag\\
 &\quad-\int_{|x'|^2<\delta^2-d^2}\frac{\alpha_2(e)}{(|x|^2+d^2)^2}\,dx'
+O\biggl(\int_{|x'|^2<\delta^2-d^2}\frac{|x'|^4}{(|x'|^2+d^2)^3}\,dx'\biggr)\notag\\
&\quad +O\biggl(\int_{|x'|^2<\delta^2-d^2}\frac{d|x'|^2}{(|x'|^2+d^2)^3}\,dx'\biggr)
+O\biggl(\int_{|x'|^2<\delta^2-d^2}\frac{|x'|^2}{(|x'|^2+d^2)^2}\,dx'\biggr)\notag\\
&=\int_{|x'|^2<\delta^2-d^2}\frac{4\alpha_2(e)x_1^2}{(|x'|^2+d^2)^3}\,dx'-
\int_{|x'|^2<\delta^2-d^2}\frac{\alpha_2(e)}{(|x|^2+d^2)^2}\,dx'+
O(d|\log d|)+O(1)
\end{align}
with
\begin{equation}
\notag
\int_{|x'|^2<\delta^2-d^2}\frac{4x_1^2}{(|x'|^2+d^2)^3}\,dx'=\frac{1}{d}\int_{\R^3}\frac{4x_1^2}{(|x'|^2+1)^3}+O(1)\,,
\end{equation}
\begin{equation}
\notag
\int_{|x'|^2<\delta^2-d^2}\frac{1}{(|x'|^2+d^2)^2}\,dx'=\frac{1}{d}\int_{\R^3}\frac{1}{(|x'|^2+1)^2}\,dx'+O(1)\,,
\end{equation}
and
\begin{equation}
\notag
\int_{\R^3}\frac{4x_1^2}{(|x'|^2+1)^3}\,dx'=\frac{4}{3}\int_{\R^3}\frac{|x'|^2}{(|x'|^2+d^2)^3}\,dx'=
\int_{\R^3}\frac{1}{(|x'|^2+1)^2}\,dx'\,.
\end{equation}
Thus, also the contribution of \eqref{4.38} remains bounded as $d\to 0$.
This completes the proof. \hfill$\Box$

\section{Critical points for the function $\tau_M$}

Let $M$ be a given $3\times 3$ real matrix. For $R\in SO(3)$,  we
define the function $\tau_M(R)=\tr(RM)$. In this section, we study
the Morse theoretical properties of this function. These results
will be applied (cf. Sections $\S5-\S8$) to the matrix $M(A_0,p_0)$
(cf. \eqref{M}), for a fixed given boundary value $A_0$ and a fixed
given point $p_0$, and are a crucial ingredient in the proofs of
Theorems 1--3 (cf., in particular \eqref{4.48}). The next two lemmas
serve to characterize the set of critical points of $\tau_M$.

\begin{lemma}
\label{L4.5} $R=R_0\in SO(3)$ is a critical point for $\tau_M$ if and only if $R_0M$ is symmetric.
\end{lemma}
\textit{Proof}: $R=R_0\in SO(3)$ is a critical point if and only if
$d\tau_M(R_0)(\xi)=\frac{d}{dt}\tau_M(\exp(t\xi)R_0)|_{t=0}=\tr(\xi
R_0M)=0$ for all $\xi\in\so(3)$. Set $B:=R_0M$ and denote by $B_+$
and $B_-$ the symmetric, and anti-symmetric parts of $B$,
respectively, i.e., $B_+=(B+B^t)/2$ and $B_-=(B-B^t)/2$. Since
$\tr(\xi B_+)=0$ for all $\xi\in\so(3)$, $d\tau_M(R_0)=0$ is
equivalent to $\tr(\xi B_-)=0$ for all $\xi\in\so(3)$. Taking
$\xi=B_-$, this implies $B_-=0$. Therefore, $R_0$ is critical if and
only if $B_-=0$, i.e., $B=R_0M$ is symmetric. \hfill$\Box$

\medskip

\begin{lemma}
\label{L4.6} Assume $\det M\ne 0$. There exists a one to one correspondence
between the set of critical points of $\tau_M$ in $SO(3)$ and the
set of symmetric $3\times 3$-matrix $B$ with $\det B=\det M$ and
$B^2=M^tM$.
\end{lemma}
\textit{Proof}: Suppose $R=R_0\in SO(3)$ is a critical point for
$\tau_M$. Set $B=R_0M$. By Lemma \ref{L4.5}, $B$ is symmetric. Moreover,
since $R_0\in SO(3)$, one has $\det B=\det M$ and
$B^2=B^tB=M^tR_0^tR_0M=M^tM$.

Conversely, suppose $B$ is a symmetric matrix with $\det B=\det M$
and $B^2=M^tM$. Define $R_0:=BM^{-1}$. Then $\det R_0=1$ and
$R_0^tR_0=1$, i.e., $R_0\in SO(3)$. Since $R_0M=B$ is symmetric,
$R_0$ is a critical point for $\tau_M$ by Lemma \ref{L4.5}. This completes
the proof. \hfill$\Box$

\medskip

In the following we assume $\det M\ne0$. See Remark \ref{R4.1} at the end
of this section for the case $\det M=0$. By Lemma \ref{L4.6}, the problem
of finding critical points of $\tau_M$ is transformed into the
problem of finding symmetric matrices $B$ with $\det B=\det M$ and
$B^2=M^tM$. The latter is easily solved as follows:

A symmetric matrix $B$ is diagonalizable by an orthogonal matrix.
So, there exists $P\in O(3)$ such that $P^{-1}BP=\begin{pmatrix} \lambda_1&0&0\\ 0&\lambda_2&0\\
0&0&\lambda_3\end{pmatrix}$, where the $\lambda_i$'s ($1\le i\le 3$) are the
eigenvalues of $B$ with $|\lambda_1|\ge|\lambda_2|\ge|\lambda_3|$.
We may assume $P\in SO(3)$ (otherwise, we take $-P$).

Denote by $\mu_1\ge\mu_2\ge \mu_3> 0$ the eigenvalues of $M^tM$
($M^tM$ is symmetric and positive). The condition $B^2=M^tM$ entails
$P\begin{pmatrix}\lambda_1^2&0&0\\ 0&\lambda_2^2&0\\
0&0&\lambda_3^2\end{pmatrix}P^{-1}=M^tM$, thus $\lambda_i^2=\mu_i$,
for $i=1,2,3$, or, equivalently,  $\lambda_i=\pm\sqrt{\mu_i}$. Since
$\det B=\lambda_1\lambda_2\lambda_3$, and $\det M = \pm
\sqrt{\mu_1\mu_2\mu_3}$, the condition $\det M= \det B$ yields the
possibilities for the matrix $B$ enlisted in the following lemma.
Notice that, in all the different cases, the critical value
corresponding to the critical point
$$R_0 := BM^{-1}=P\begin{pmatrix} \lambda_1&0&0\\ 0&\lambda_2&0\\
0&0&\lambda_3\end{pmatrix}P^{-1} M^{-1}$$ is given by $\tr \,
B=\lambda_1+\lambda_2+\lambda_3$.

\begin{lemma}
\label{L4.7} Let $M$ be a given $3\times 3$ real matrix, such that  $M^tM$ has eigenvalues $\mu_1\ge\mu_2\ge \mu_3>
0$, and let $B$ a symmetric $3\times 3$ real matrix, with $\det
B=\det M$ and $B^2=M^tM$. Then, there exists $P\in SO(3)$ such that
$B=P\begin{pmatrix} \lambda_1&0&0\\ 0&\lambda_2&0\\
0&0&\lambda_3\end{pmatrix}P^{-1}$, where the possible $\lambda_i$'s
are enlisted below:

\begin{enumerate}[(1)]
\item If $\det M:= +\sqrt{\mu_1\mu_2\mu_3}>0$. In this case $\lambda_1\lambda_2\lambda_3>0$, and there are four possibilities:
\begin{enumerate}
\item $\lambda_1=\sqrt{\mu_1}$, $\lambda_2=\sqrt{\mu_2}$,
$\lambda_3=\sqrt{\mu_3}$, yielding the critical value
$\sqrt{\mu_1}+\sqrt{\mu_2}+\sqrt{\mu_3}$,
\item $\lambda_1=\sqrt{\mu_1}$, $\lambda_2=-\sqrt{\mu_2}$, $\lambda_3=-\sqrt{\mu_3}$, yielding the critical value
$\sqrt{\mu_1}-\sqrt{\mu_2}-\sqrt{\mu_3}$,
\item $\lambda_1=-\sqrt{\mu_1}$, $\lambda_2=\sqrt{\mu_2}$,
$\lambda_3=-\sqrt{\mu_3}$, yielding the  critical value
$-\sqrt{\mu_1}+\sqrt{\mu_2}-\sqrt{\mu_3}$,
\item $\lambda_1=-\sqrt{\mu_1}$, $\lambda_2=-\sqrt{\mu_2}$,
$\lambda_3=\sqrt{\mu_3}$, yielding the  critical value
$-\sqrt{\mu_1}-\sqrt{\mu_2}+\sqrt{\mu_3}$.
\end{enumerate}
\item If $\det M:=-\sqrt{\mu_1\mu_2\mu_3}<0$. In this case $\lambda_1\lambda_2\lambda_3<0$, and there also are four possibilities:
\begin{enumerate}
\item $\lambda_1=\sqrt{\mu_1}$, $\lambda_2=\sqrt{\mu_2}$, $\lambda_3=-\sqrt{\mu_3}$,
yielding the critical value
$\sqrt{\mu_1}+\sqrt{\mu_2}-\sqrt{\mu_3}$,
\item $\lambda_1=\sqrt{\mu_1}$, $\lambda_2=-\sqrt{\mu_2}$,
$\lambda_3=\sqrt{\mu_3}$, yielding the critical value
$\sqrt{\mu_1}-\sqrt{\mu_2}+\sqrt{\mu_3}$,
\item $\lambda_1=-\sqrt{\mu_1}$, $\lambda_2=\sqrt{\mu_2}$,
$\lambda_3=\sqrt{\mu_3}$, yielding the critical value
$-\sqrt{\mu_1}+\sqrt{\mu_2}+\sqrt{\mu_3}$,
\item $\lambda_1=-\sqrt{\mu_1}$, $\lambda_2=-\sqrt{\mu_2}$,
$\lambda_3=-\sqrt{\mu_3}$, yielding the critical value
$-\sqrt{\mu_1}-\sqrt{\mu_2}-\sqrt{\mu_3}$.
\end{enumerate}
\end{enumerate}
\end{lemma}

\medskip

(Note that the list above gives all the possible critical values for
$\tau_M$, if $\det M\neq 0.)$

We next study the non-degeneracy of all the critical points and
calculate their Morse indexes, still in the hypothesis $\det M\neq
0.$

To this purpose, we study the  Hessian of $\tau_M$ at $R_0:= B
M^{-1}$. This  is given by
$$\nabla^2\tau_M(R_0)(\xi,\zeta)=\frac{\partial^2}{\partial s\partial t}\tau_M(\exp(s\xi+t\zeta)R_0M)\big|_{s=t=0}=
\frac{1}{2}\tr((\xi\zeta+\zeta\xi)R_0M)=\tr(\xi\zeta B)$$ for
$\xi,\zeta\in\so(3)$, where in the last equality we have used the
fact that $R_0M$ is symmetric.

Let $P\in SO(3)$ be as before.  We represent $\nabla^2\tau_M(R_0)$
with respect to the basis $P\xi_iP^{-1}$
($i=1,2,3$), with $\xi_1=\begin{pmatrix} 0&0&0\\
0&0&-1\\ 0&1&0\end{pmatrix}$, $\xi_2=\begin{pmatrix} 0&0&1\\ 0&0&0\\
-1&0&0\end{pmatrix}$, $\xi_3=\begin{pmatrix} 0&-1&0\\ 1&0&0\\
0&0&0\end{pmatrix}$.

\noindent An easy calculation shows that
\begin{align*}
\nabla^2\tau_M(R_0)(P\xi_iP^{-1},P\xi_jP^{-1})=\tr(\xi_i\xi_jP^{-1}BP)=0,\qquad 1\le i\ne j\le 3\;,\\
\nabla^2\tau_M(R_0)(P\xi_1P^{-1},P\xi_1P^{-1})=\tr(\xi_1^2P^{-1}BP)=-\lambda_2-\lambda_3\;,\\
\nabla^2\tau_M(R_0)(P\xi_2P^{-1},P\xi_2P^{-1})=\tr(\xi_2^2P^{-1}BP)=-\lambda_1-\lambda_3\;,\\
\nabla^2\tau(R_0)(P\xi_3P^{-1},P\xi_3P^{-1})=\tr(\xi_3^2P^{-1}BP)=-\lambda_1-\lambda_2\;.\notag\\
\end{align*}
So, the Hessian $\nabla^2\tau_M(R_0)$ is written as
$\begin{pmatrix} -\lambda_2-\lambda_3&0&0\\
0&-\lambda_1-\lambda_3&0\\ 0&0&-\lambda_1-\lambda_2\end{pmatrix}.$

From this, the Morse indexes of the critical points corresponding to
all the possible cases, which we  described in Lemma \ref{L4.7}, are easily
computed as stated in the following lemma.

\begin{lemma}
\label{L4.8} Let $M$, $B$, $\tau_M$ be given as in the previous
lemmas.
\begin{enumerate}[(1)]
\item If $\det M:=+\sqrt{\mu_1\mu_2\mu_3}>0$. In the four different cases, (1)-(a), (b), (c), (d)
of Lemma \ref{L4.7}, the following assertions hold, respectively:
\begin{enumerate}
\item all the $\lambda_j$'s are strictly positive, the Hessian is negative-definite, the critical points are non-degenerate
with Morse index equal to $3$,
\item $-\lambda_1-\lambda_2=-\sqrt{\mu_1}+\sqrt{\mu_2}\le 0$,
$-\lambda_2-\lambda_3=\sqrt{\mu_2}+\sqrt{\mu_3}>0$,
$-\lambda_1-\lambda_3=-\sqrt{\mu_1}+\sqrt{\mu_3}\le 0$, thus the
critical points are non-degenerate if and only if $\mu_1>\mu_2$, in
which case, they have Morse index equal to $2$,
\item $-\lambda_1-\lambda_2=\sqrt{\mu_1}-\sqrt{\mu_2}\ge 0$, $-\lambda_2-\lambda_3=-\sqrt{\mu_2}+\sqrt{\mu_3}\le 0$,
$-\lambda_1-\lambda_3=\sqrt{\mu_1}+\sqrt{\mu_3}>0$, thus the
critical points are non-degenerate if and only if
$\mu_1>\mu_2>\mu_3$, in which case, they have Morse index equal to
$1$,
\item $-\lambda_1-\lambda_2=\sqrt{\mu_1}+\sqrt{\mu_2}>0$, $-\lambda_2-\lambda_3=\sqrt{\mu_2}-\sqrt{\mu_3}\ge 0$,
$-\lambda_1-\lambda_3=\sqrt{\mu_1}-\sqrt{\mu_3}\ge 0$, thus the
critical points are non-degenerate if and only if $\mu_2>\mu_3$, in
which case, they have Morse index equal to $0$;
\end{enumerate}
\item if $\det M:=-\sqrt{\mu_1\mu_2\mu_3}<0$. In the four different cases (2)-(a), (b), (c), (d) of Lemma \ref{L4.7}, the following assertions hold, respectively:
\begin{enumerate}
\item $-\lambda_1-\lambda_2=-\sqrt{\mu_1}-\sqrt{\mu_2}<0$, $-\lambda_2-\lambda_3=-\sqrt{\mu_2}+\sqrt{\mu_3}\le 0$,
$-\lambda_1-\lambda_3=-\sqrt{\mu_1}+\sqrt{\mu_3}\le 0$, thus the
critical points are non-degenerate if and only if $\mu_2>\mu_3$, in
which case, they have Morse index equal to $3$,
\item  $-\lambda_1-\lambda_2=-\sqrt{\mu_1}+\sqrt{\mu_2}\le 0$, $-\lambda_2-\lambda_3=\sqrt{\mu_2}-\sqrt{\mu_3}\ge 0$,
$-\lambda_1-\lambda_3=-\sqrt{\mu_1}-\sqrt{\mu_3}< 0$, thus  the
critical points are non-degenerate if and only if
$\mu_1>\mu_2>\mu_3$, in which case, they have  Morse index equal to
$2$,
\item $-\lambda_1-\lambda_2=\sqrt{\mu_1}-\sqrt{\mu_2}\ge 0$, $-\lambda_2-\lambda_3=-\sqrt{\mu_2}-\sqrt{\mu_3}< 0$,
$-\lambda_1-\lambda_3=\sqrt{\mu_1}-\sqrt{\mu_3}\ge 0$, thus the
critical points are non-degenerate if and only if $\mu_1>\mu_2$, in
which case, they have Morse index equal to $1$,
\item all $\lambda_i$'s are strictly negative, and the Hessian is positive-definite, thus the corresponding critical points
are non-degenerate, with Morse index equal to $0$.
\end{enumerate}
\end{enumerate}
\end{lemma}

\medskip

We finally prove the following lemma.
\begin{lemma}
\label{L4.9} For all the
non-degenerate cases (described in Lemma \ref{L4.8}), there corresponds
exactly one critical point for each critical value.
\end{lemma}
\textit{Proof}: Note that in the non-degenerate cases $\lambda_i +
\lambda_j\neq 0$, $\forall i, j$.  We must show that for a given
critical value (as listed in Lemma \ref{L4.7}), there exists exactly one
$B$ satisfying $\det B=\det M$ and $B^2=M^tM$. So suppose both of
$B$ and $B'$ satisfy these two conditions. Then, there exist $P,Q\in
SO(3)$ such that
$P^{-1}BP=D$ and $Q^{-1}B'Q=D'$, where $D=\begin{pmatrix}\lambda_1&0&0\\ 0&\lambda_2&0\\ 0&0&\lambda_3\end{pmatrix}$ and
$D'=\begin{pmatrix}\lambda_1'&0&0\\
0&\lambda_2'&0\\ 0&0&\lambda_3'\end{pmatrix}$ are diagonal matrices
with $|\lambda_1|\ge|\lambda_2|\ge|\lambda_3|$, and
$|\lambda_1'|\ge|\lambda_2'|\ge|\lambda_3'|$. By Lemma \ref{L4.7}, $D=D'$.
The condition $B^2={B'}^2=M^tM$ implies
$PD^2P^{-1}=Q{D'}^2Q^{-1}=QD^2Q^{-1}$, thus $CD^2C^{-1}=D^2$, with
$C:=Q^{-1}P\in SO(3)$. Denoting by $c_{ij}$ the entries of $C$, we
thus obtain
$$(\lambda_i^2-\lambda_j^2)c_{ij}=0\;, \quad 1\le i,j\le 3\;.$$
Since $\lambda_i+\lambda_j\ne 0$, for all $i,j$, in all the
non-degenerate cases, this yields
$$(\lambda_i-\lambda_j)c_{ij}=0\,,\quad  1\le i,j\le 3\;.$$
But this last condition is equivalent to $CDC^{-1}=D$, thus  $B=B'$,
and the lemma is proved. \hfill$\Box$
\medskip

\newtheorem{remark}{Remark}[section]
\begin{remark}
\label{R4.1}
If $\det M=0$, similar arguments show that the possible critical
values for $\tau_M$ are: (a) $\sqrt{\mu_1}+\sqrt{\mu_2}$, (b)
$\sqrt{\mu_1}-\sqrt{\mu_2}$, (c) $-\sqrt{\mu_1}+\sqrt{\mu_2}$, and
(d) $-\sqrt{\mu_1}-\sqrt{\mu_2}$. Moreover, for each critical value
there corresponds exactly one critical point in the cases (a), (d),
and, also, in the cases (b),(c), provided that $\sqrt{\mu_1}\neq
\sqrt{\mu_2}$. The corresponding critical point is non-degenerate if
and only if $\mu_1\geq\mu_2>0$ in the cases (a),(d), and if and only
if $\mu_1>\mu_2>0$ in the cases (b),(c). The corresponding Morse
indexes are $3$ in the case (a), $2$ in the case (b), $1$ in the
case (c), $0$ in the case (d). Thus, $\tau_M$ is a Morse function
exactly when $\mu_1>\mu_2>0$. In this case (and this is the only one
we need in this paper), the above four values are in fact critical
values for $\tau_M$. Indeed, by the Ljusternik-Schnirelmann theory
(see~\cite{Struwe} and the proof of Theorem \ref{3}), any
function on $SO(3)$ has at least four critical points, since the
Ljusternik-Schnirelmann category of $SO(3)\cong \R P^3$ is $4$
(see~\cite{Struwe}).
\end{remark}

\section{Asymptotic estimates of $\Je(\q)$ and $\Je'(\q)$}

In order to prove Theorems 1--3, we need the following lemmas which compare asymptotically, as $\epsilon\to 0$,
the functional
$\Je(\q)=\epsilon^2\YMe(A(\q)+ a(\q))$ with the functional
$J_{\epsilon}(\q)=\epsilon^2\YMe(A(\q))$,
both defined on the parameter space
$\PP(d_0,\lambda_0;D_1,D_2;\epsilon)$ (cf.
\eqref{curl}-\eqref{par}). Since an estimate of $J_{\epsilon}(\q)$
and of its derivative $J^\prime_{\epsilon}(\q)$ are given in $\S3.2$
of  \cite{IM1}, these lemmas yield estimates for $\Je(\q)$ and its
derivative $\Je^\prime(\q).$

\medskip\noindent
We recall that for connections $A$ on the bundle $P$
and one-forms $a\in
C^{\infty}(T^{\ast}\overline{B}^4\otimes\Ad(P))$, the
$L^2_1$-norm $\|a\|_{A;1,2}$ is defined by
\begin{equation}
\label{norm} \|a\|_{A;1,2}:=\|{\nabla_A}^{\epsilon}a\|_2+\|a\|_2\;,
\end{equation}
where
$\|\cdot\|_2$ is
the $L^2$-norm on $B^4$. Observe that the space $L^2_{0,1}(T^{\ast}B^4\otimes\Ad(P))$ (the completion of
$C^{\infty}_0(T^{\ast}B^4\otimes\Ad(P))$ with respect to the norm
above) is independent of the choice of
the connection $A$. We also recall that the \emph{dual} $L^2_1$-norm of $\nabla\YMe$ is defined by
\begin{equation}
\label{dualnorm}
\|\nabla\YMe(A)\|_{A;1,2,\ast}:=\sup\{\nabla\YMe(A)(a):a\in L^2_{0,1}(T^{\ast}B^4\otimes\Ad(P)),~\|a\|_{A;1,2}\le 1\}\;.
\end{equation}

\begin{lemma}
\label{L4.1}For $\q\in\PP(d_0,\lambda_0;D_1,D_2;\epsilon)$, there holds
\begin{equation}
\label{4.1}
\Je(\q)=J_{\epsilon}(\q)+r_3(\q),
\end{equation}
where $|r_3(\q)|\lesssim\epsilon^3$ uniformly with respect to $\q$.
\end{lemma}
\textit{Proof:}
We have
\begin{align}
\label{4.2}
\Je(\q)&=\epsilon^2\int_{B^4}|{F_{A(\q)+a(\q)}}^{\epsilon}|^2\,dx\notag\\
&=\epsilon^2\int_{B^4}|{F_{A(\q)}}^{\epsilon}|^2\,dx+\epsilon^2\int_{B^4}|{d_{A(\q)}}^{\epsilon}a(\q)|^2\,dx+
\frac{\epsilon^4}{4}\int_{B^4}|[a(\q),a(\q)]|^2\,dx\notag\\
+2\epsilon^2&\int_{B^4}({F_{A(\q)}}^{\epsilon},{d_{A(\q)}}^{\epsilon}a(\q))\,dx+\epsilon^3\int_{B^4}({F_{A(\q)}}^{\epsilon},[a(\q),a(\q)])\,dx
+\epsilon^3\int_{B^4}({d_{A(\q)}}^{\epsilon}a(\q),[a(\q),a(\q)])\,dx.
\end{align}
The first term on the right hand side of \eqref{4.2} is
$J_{\epsilon}(\q)$. The remaining terms are easily estimated by
Lemmas 3.2 and 3.9 in \cite{IM1}, and the Sobolev inequality: the second, third, fifth
and sixth terms are bounded by $C\epsilon^3$ for some $C>0$ depending
only on $d_0,\lambda_0,D_1$ and $D_2$. The fourth term is estimated
as
$$\epsilon^2\biggl|\int_{B^4}({F_{A(\q)}}^{\epsilon},{d_{A(\q)}}^{\epsilon}a(\q))\,dx\biggr|\le\epsilon^2\|\nabla\YMe(A(\q))\|_{A(\q);1,2,\ast}\|a(\q)\|_{A(\q);1,2;B^4}
\lesssim\epsilon^3.$$
Combining all these estimates, \eqref{4.1} follows
easily. \hfill$\Box$

\medskip

The following lemma compares the derivative of
$\Je(\q)$ with the derivative of $J_{\epsilon}(\q)$.
We use the following notation: $\q_i(\q)$ for $i=1,..., 8$ are the vector fields constructed in \cite{IM3} such that
$a_i(\q) = A_{\q_i}(\q)$ (the directional derivative of $A(\q)$ in the direction of $\q_i(\q)).$

\begin{lemma}
\label{L4.2} The following holds:
\begin{equation}
\notag \langle\Je'(\q),\q_i(\q)\rangle=\langle
J_{\epsilon}'(\q),\q_i(\q)\rangle+r_{4,i}(\q),
\end{equation}
where $|r_{4,i}(\q)|\lesssim\epsilon^4$ for $1\le i\le 4$ and
$|r_{4,i}(\q)|\lesssim\epsilon^{7/2}$ for $5\le i\le 8$ uniformly
with respect to $\q\in\PP(d_0,\lambda_0;D_1,D_2;\epsilon)$.
\end{lemma}
\textit{Proof:} One has
\begin{align}
\label{4.4}
\langle\mathcal{J}_{\epsilon}'(\q),\q_i(\q)\rangle&=\epsilon^2\langle\nabla\YMe(A(\q)+a(\q)),\A_i(\q)+
a_{\q_i}(\q)\rangle\notag\\
&=\epsilon^2\langle\nabla\YMe(A(\q)),\A_i(\q)\rangle+\epsilon^2\langle\nabla\YMe(A(\q)),a_{\q_i}(\q)\rangle\notag\\
&\quad+\epsilon^2\langle\nabla^2\YMe(A(\q))a(\q),\A_i(\q)\rangle+
\epsilon^2\langle\nabla^2\YMe(A(\q))a(\q),a_{\q_i}(\q)\rangle\notag\\
&\quad+\epsilon^2\langle R(\q;a(\q)),\A_i(\q)\rangle+\epsilon^2\langle
R(\q;a(\q)),a_{\q_i}(\q)\rangle.
\end{align}
We estimate each terms in \eqref{4.4}. The first term is  $\langle
J'_{\epsilon}(\q),\q_i(\q)\rangle$. By Lemmas 3.2 and 3.10 in \cite{IM1}, the
second term is estimated as
\begin{align}
|\epsilon^2\langle\nabla\YMe(A(\q)),a_{\q_i}(\q)
\rangle|&\le\epsilon^2\|\nabla\YMe(A(\q))\|_{A(\q);1,2,\ast}\|a_{\q_i}(\q)\|_{A(\q);1,2;B^4}\notag\\
&\lesssim\epsilon^4\quad\text{for $1\le i\le 4$}; \quad\text{or $\quad\lesssim\epsilon^{7/2}$}\quad
\text{for $5\le i\le 8$}.
\end{align}
By Lemmas 3.2,  3.7, 3.8, 3.9 in  \cite{IM1}, the third term is
estimated as
\begin{align}
&|\epsilon^2\langle\nabla^2\YMe(A(\q))a(\q),\A_i(\q)
\rangle|\le\epsilon^2\|\nabla^2\YMe(A(\q))\A_i(\q)\|_{A(\q);1,2,\ast}\|a(\q)\|_{A(\q);1,2;B^4}\notag\\
&\lesssim\epsilon^{5/2}\|(\nabla^2\YMe(A(\q))-\nabla^2\YMe(\tilde{A}(\q)))\A_i(\q)\|_{A(\q);1,2,\ast}
+\epsilon^{5/2}\|\nabla^2\YMe(\tilde{A}(\q))(\A_i(\q)-\tilde{\A}_i(\q))\|_{A(\q);1,2,\ast}\notag\\
&\lesssim\epsilon^4+\epsilon^{5/2}\|\A_i(\q)-\tilde{\A}_i(\q)\|_{A)(\q);1,2;B^4}
\lesssim\epsilon^{4}\quad\text{ for } 1\le i\le 4\;; \quad\text{or}\;
\lesssim\epsilon^{7/2}\quad \text{ for } 5\le i\le 8\,.
\end{align}
By Lemmas 3.2, 3.9 and 3.10 in \cite{IM1}, the fourth term is estimated as
\begin{align}
&|\epsilon^2\langle\nabla^2
\YMe(A(\q))a(\q),a_{\q_i}(\q)\rangle|\le\epsilon^2\|a(\q)\|_{A(\q);1,2;B^4}\|a_{\q_i}(\q)\|_{A(\q);1,2;B^4}\notag\\
&\lesssim\epsilon^{4}\quad\text{ for } 1\le i\le 4\;; \quad\text{or}\;
\lesssim\epsilon^{7/2} \quad \text{ for } 5\le i\le 8\;,
\end{align}
Similarly, by Lemmas 3.2, 3.3, 3.9 in \cite{IM1}, the fifth term is estimated as
\begin{align}
|\epsilon^2\langle R(\q;a(\q)),\A_i(\q)\rangle|&\lesssim\epsilon^2(\epsilon\|a(\q)\|^2_{A(\q);1,2;B^4}+
\epsilon^2\|a(\q)\|^3_{A(\q);1,2;B^4})\notag\\
&\lesssim\epsilon^4\quad \text{for } 1\le i\le 8\,,
\end{align}
and, by Lemmas 3.2, 3.3, 3.9, 3.10 in \cite{IM1}, the last term is estimated as
\begin{align}
\label{4.9}
|\epsilon^2\langle R(\q;a(\q)),a_{\q_i}(\q)\rangle|&\lesssim\epsilon^2(\epsilon\|a(\q)\|^2_{A(\q);1,2;B^4}+
\epsilon^2\|a(\q)\|^3_{A(\q);1,2;B^4})\|a_{\q_i}(\q)\|_{A(\q);1,2;B^4}\notag\\
&\lesssim\epsilon^5\quad \text{for } 1\le i\le 8\,.
\end{align}
The Lemma follows from \eqref{4.4}--\eqref{4.9}. \hfill$\Box$

\medskip

\noindent By Proposition 3.1 in \cite{IM1} and Lemma 3.2 in \cite{IM3}, the leading term $\langle
J'_{\epsilon}(\q),\q_i(\q)\rangle$ of equation \eqref{4.1} is estimated
as
\begin{equation}
\notag
\langle J'_{\epsilon}(\q),\q_1(\q)\rangle =a_{11}(\q)\Big\langle
J'_{\epsilon}(\q),\frac{\partial}{\partial p_1}\Big\rangle
\simeq\epsilon^{7/2}C_1(\q)\,,
\end{equation}
where $C_1(\q)$ is a constant depending only on $\q$.

\noindent Similarly, it follows
that:
\begin{equation}
\notag
\langle
J'_{\epsilon}(\q),\q_i(\q)\rangle\simeq\epsilon^{7/2}C_i(\q)\quad
\text{for} \quad 1\le i\le 4,
\end{equation}
and
\begin{equation}
\notag
\langle
J'_{\epsilon}(\q),\q_i(\q)\rangle\simeq\epsilon^3C_i(\q)\quad \text{
for} \quad 5\le i\le 8,
\end{equation}
with $C_i(\q)$ depending only on $\q$.

\section{Proof of Theorem 1}

We are now ready to prove the first of our existence theorems stated
in \cite{IM1}, which we state again here for the convenience of the
reader. We recall that Theorems 1-3 can all be restated in terms of
the Dirichlet problem for the $SU(2)$-Yang Mills functional with
boundary value $\epsilon A_0$ (cf. $\S$1 of this paper or $\S$2.2 of
\cite{IM1}).
\newtheorem{theorem}{Theorem}
\begin{theorem}
\label{1}
Let us define the function $G_1^\pm (p):=\frac{(\sqrt{\mu_1(p)}+ \sqrt{\mu_2(p)}\pm\sqrt{\mu_3(p)})^2}{F(p)}\;,\; p\in B^4,$
and assume that $p_0\in B^4$ satisfies either of the following hypotheses (1),(2):
\begin{enumerate}[(1)]
\item $\det M(A_0,p_0)>0$  and $p_0$ is an isolated local maximum point of $G_1^+(p);$
\item $\det M(A_0,p_0)<0$ and $p_0$ is an isolated local maximum point of $G_1^-(p).$
\end{enumerate}
Then, there exists $\epsilon_0>0$ and  a family of connections
$\{{A}_\epsilon\}$ indexed by $\epsilon\in (0, \epsilon_0]$ with the
following properties: ${A}_{\epsilon}$ is a solution to
$\bigl(\mathcal{D}_{\epsilon}\bigr)$ in $\mathcal{A}_{+1}(A_0)$;
$\epsilon^2|{F_{{A}_{\epsilon}}}^{\epsilon}|^2\,dx\to
8\pi^2\delta_{p_0}$ as $\epsilon\to 0$ in the sense of measures
(i.e. $\epsilon {F_{{A}_{\epsilon}}}^{\epsilon}$ concentrates at
$p_0$ as $\epsilon\to 0$).
\end{theorem}

\medskip
\noindent\textit{Proof:} \emph{Case (1).} We assume that $p_0\in B^4$ is such that
$\det M(A_0,p_0)>0$, and is an isolated local maximum point for the
function $G_1^+$. Thus, there exists $\delta>0$ such that
$B_{\delta}(p_0)\Subset B^4$, $\det M(A_0,p)>0$ for all $p\in
B_{\delta}(p_0)$, and $G_1^+(p_0)>G_1(p)^+$ for all $p\in
B_{\delta}(p_0)\setminus\{p_0\}$. From now on, we fix such $\delta$
and restrict our choice of $0<D_1<D_2$ as follows:
\begin{align}
\label{4.42}
0<D_1<\frac{1}{2}\frac{\Gamma_1^+(p_0)}{F(p_0)},\\
\label{4.43}
D_1^2F(p)-2D_1\Gamma_1^+(p)\ge-\frac{1}{2}G_1^+(p_0),\\
\label{4.44}
2\frac{\Gamma_1^+(p_0)}{F(p_0)}<D_2,\\
\label{4.45}
D_2^2F(p)-2D_2\Gamma_1^+(p)\ge-\frac{1}{2}G_1^+(p_0),
\end{align}
for all $p\in B_{\delta}(p_0)$, where
$\Gamma_1^+(p):=\sqrt{\mu_1(A_0,p)}+\sqrt{\mu_2(A_0,p)}+\sqrt{\mu_3(A_0,p)}.$
Note that \eqref{4.42}, \eqref{4.43} are both satisfied if $D_1>0$ is chosen
suitably small, while \eqref{4.44}, \eqref{4.45} are both satisfied if $D_2$
is chosen suitably large, since $F(p)>0$ in $B^4$ (cf. (1) of Lemma
\ref{L4.3}).

We then choose $d_0>0$ such that $B_{\delta}(p_0)\Subset
B_{1-d_0}(0)$. Under the assumptions above on $D_1$, $D_2$, $d_0$, we
define $\PP(p_0,\delta;d_0, \lambda_0;D_1,D_2;\epsilon)$, a subset
of $\PP(d_0,\lambda_0;D_1,D_2;\epsilon)$, as follows:
\begin{equation}
\notag
\PP(p_0,\delta;d_0,\lambda_0;D_1,D_2;\epsilon)=\{(p,[g],\lambda)\in\PP(d_0,\lambda_0;D_1,D_2;\epsilon):p\in
B_{\delta}(p_0)\}.
\end{equation}
Since $\q$ is a critical point for $\mathcal{J}_{\epsilon}(\q)$ if
and only if $A(\q) + a(\q)$ is a Yang Mills connection (cf.
Proposition 3.2 in \cite{IM1}), we look for critical points of the
function $\mathcal{J}_{\epsilon}$ in the interior of
$\PP(p_0,\delta;d_0,\lambda_0;D_1,D_2;\epsilon)$. Since the closure
of $\PP(p_0,\delta;d_0,\lambda_0;D_1,D_2;\epsilon)$ in
$\PP(d_0,\lambda_0;D_1,D_2;\epsilon)$ is compact, there exists a
value of the parameter
$\q_m\in\overline{\PP(p_0,\delta;d_0,\lambda_0;D_1,D_2;\epsilon)}$
such that $\mathcal{J}_{\epsilon}$ attains its minimum in
$\overline{\PP(p_0,\delta;d_0,\lambda_0;D_1,D_2;\epsilon)}$ at
$\q_m$. We show that
$\q_m\in\PP(p_0,\delta;d_0,\lambda_0;D_1,D_2;\epsilon)$ and that it
is a critical point for $\mathcal{J}_{\epsilon}$  in
$\PP(d_0,\lambda_0;D_1,D_2;\epsilon)$.

To see this, set $\q_0=\Big(\sqrt{\frac{\Gamma_1^+
(p_0)}{F(p_0)}\epsilon},[g_0],p_0\Big)\in\PP(p_0,\delta;d_0,\lambda_0;D_1,D_2;\epsilon)$,
where $[g_0]\in SO(3)$ is the maximum point of the function
\begin{equation}
\label{4.47}
SO(3)\ni[g]\mapsto\int_{B^4}\bigl((d\underline{A}_0)^-,g(dh_{p_0})^-g^{-1}\bigr)\,dx\,.
\end{equation}
Making the identification $\text{Im}\,\HH\cong\R^3$, the function
\eqref{4.47} can be rewritten as
\begin{align}
\label{4.48} SO(3)\ni
R\mapsto&\int_{B^4}\bigl((d\underline{A}_0)^-,R(dh_{p_0})^-\bigr)\,dx
=\sum_{1\le i,j\le 3}R_{ij}\int_{B^4}\bigl((d\underline{A}_{0,
i})^-,(dh_{p_0,j})^-\bigr)\,dx =\text{tr}(RM(A_0,p_0))\,,
\end{align}
where $M(A_0,p)$ is the matrix defined in \eqref{M}.

\noindent The results obtained in $\S3$ for the case (1)-(a)
yield
\begin{equation}
\label{4.49} \int_{B^4}\bigl((d\underline{A}_{0,
i})^-,g_0(dh_p)^-g_0^-\bigr)\,dx =\Gamma_1^+(p).
\end{equation}
By Lemma \ref{L4.1} and by the asymptotic expansion in Proposition
3.1 in \cite{IM1}, one obtains
\begin{equation}
\label{4.50}
\Je(\q_0)=C_{\epsilon}-2\frac{\Gamma_1^+(p_0)^2}{F(p_0)}\epsilon^2+r_5(\q)
=C_{\epsilon}-2G_1^+(p_0)\epsilon^2+r_{5,\epsilon},\quad\mbox{ with } |r_{5,\epsilon}|\lesssim\epsilon^3\,,
\end{equation}
where
$C_{\epsilon}=8\pi^2+\epsilon^2\int_{B^4}|{F_{\underline{A}_{\epsilon}}}^{\epsilon}|^2\,dx$
is a constant depending only on $\epsilon$.
On the other hand, writing $\q_m=(p_m,[g_m],\lambda_m)$,
the results obtained in $\S3$ for the case (1)-(a) yield
\begin{align}
\label{4.51}
&\Je(\q_m)=C_{\epsilon}+2\lambda_m^4F(p_m)-4\epsilon\lambda_m^2
\int_{B^4}\bigl((d\underline{A}_0)^-,g_m(dh_{p_m})^-g_m^{-1}\bigr)\,dx +r_{6,\epsilon}\notag\\
&\ge
C_{\epsilon}+2\lambda_m^4F(p_m)-4\epsilon\lambda_m^2\Gamma_1^+(p_m)+r_{6,\epsilon}
\ge C_{\epsilon}-2\epsilon^2G_1^+(p_m)+r_{6,\epsilon}, \quad\mbox {
with }|r_{6,\epsilon}|\lesssim\epsilon^3\,,
\end{align}
where we have used the
fact that the function
$\lambda\mapsto2\lambda^4F(p_m)-4\epsilon\lambda^2\Gamma_1^+(p_m)$ has
minimum value $-2\epsilon^2G_1^+(p_m)$.

\noindent Since $\Je(\q_0)\ge\Je(\q_m)$,  from \eqref{4.50}, \eqref{4.51} we
derive
\begin{equation}
\label{4.52}
G_1^+(p_m)\ge
G_1^+(p_0)+\epsilon^{-2}(r_{6,\epsilon}-r_{5,\epsilon}).
\end{equation}
Setting $\gamma:=\min_{p\in\partial
B_{\delta}(p_0)}(G_1^+(p_0)-G_1^+(p))$, by \eqref{4.52} we would have
\begin{equation}
\notag
0<\gamma\le \epsilon^{-2}(r_{5,\epsilon}-r_{6,\epsilon})
\end{equation}
if $p_m\in\partial B_{\delta}(p_0)$. But this is a contradiction for
small $\epsilon>0$, since
$|\epsilon^{-2}(r_{6,\epsilon}-r_{5,\epsilon})|\lesssim\epsilon$. Thus, $p_m\in B_{\delta}(p_0)$.

To prove $D_1\epsilon<\lambda_m<D_2\epsilon$, suppose first that $\lambda_m^2=D_1\epsilon$. In this case,
 by \eqref{4.43}, \eqref{4.50}, one has,  for small positive $\epsilon$,
 \begin{align}
 \notag
 \Je(\q_m)&=C_{\epsilon}+2\lambda_m^4F(p_m)-
 4\epsilon\lambda_m^2\int_{B^4}\bigl((d\underline{A}_0)^-,g_m(dh_{p_m})^-g_m^{-1}\bigr)\,dx +r_{6,\epsilon}\notag\\
 &\ge C_{\epsilon}+2\lambda_m^4F(p_m)-4\epsilon\lambda_m^2\Gamma_1^+(p_m)+r_{6,\epsilon}\notag\\
&\ge C_{\epsilon}+2\epsilon^2(D_1^2F(p_m)-2D_1\Gamma_1^+(p_m))+r_{6,\epsilon}\ge C_{\epsilon}-\epsilon^2G_1^+(p_0)+r_{6,\epsilon}>\Je(\q_0)\,.
 \end{align}

This contradicts the minimality of $\q_m$. Thus, we have
$D_1\epsilon<\lambda_m^2$. A similar argument shows that
 $\lambda_m^2<D_2\epsilon$, for small $\epsilon>0$.

Summing up, for small $\epsilon>0$, one has that
$\q_m\in\PP(p_0,\delta;d_0,\lambda_0;D_1,D_2;\epsilon)$ and
is a critical point of $\Je$ in
$\PP(d_0,\lambda_0;D_1,D_2;\epsilon)$. Therefore, by  Proposition
3.2 in \cite{IM1}, $A_{\epsilon}:=A(\q_m)+a(\q_m)$ is a solution to the
Dirichlet problem  $(\mathcal{D}_{\epsilon})$.
Moreover, our construction yields
$\epsilon^2|{F_{A_{\epsilon}}}^{\epsilon}|^2\,dx\to 8\pi^2\delta_{p_0}$ in the
sense of Radon measures, as $\epsilon\to 0$. This completes the
proof of case (1) of Theorem 1.

\noindent \emph{Case (2).} The proof of case (2) is quite similar to the one just given.
The only difference is that the maximum value for the function
\eqref{4.47} is
$\Gamma^-_1(p):=\sqrt{\mu_1(A_0,p)}+\sqrt{\mu_2(A_0,p)}-\sqrt{\mu_3(A_0,p)}$ (cf.
(2)-(a) in $\S3$). \hfill$\Box$

\section{Proof of Theorem \ref{2}}

In this section we prove the second of our existence theorems, i.e., the following

\begin{theorem}
\label{2}
Let us define the functions $G_2^\pm(p):=\frac{(\sqrt{\mu_1(p)}-\sqrt{\mu_2(p)}\mp\sqrt{\mu_3(p)})^2}{F(p)};$
\hfill

\noindent $G_3^-(p):=\frac{(-\sqrt{\mu_1(p)}+\sqrt{\mu_2(p)}+\sqrt{\mu_3(p)})^2}{F(p)};$ $G_1^0(p):=\frac{(\sqrt{\mu_1(p)}+\sqrt{\mu_2(p)})^2}{F(p)};$ $G_2^0(p):=\frac{(\sqrt{\mu_1(p)}-\sqrt{\mu_2(p)})^2}{F(p)}.$

\noindent
Assume that $p_0\in B^4$ satisfies one of the following conditions (1)-(a),(b), (2)-(a),(b),(c), (3)-(a),(b):
\begin{enumerate}[(1)]
\item $\det M(A_0,p_0)>0$ and
\begin{enumerate}
\item $p_0$ is a non-degenerate critical point of $G_1^+(p)$, or
\item $\sqrt{\mu_1(A_0,p_0)}>\sqrt{\mu_2(A_0,p_0)}+\sqrt{\mu_3(A_0,p_0)}$ and $p_0$ is a non-degenerate critical point of $G_2^+(p);$
\end{enumerate}
\item $\det M(A_0,p_0)<0$ and
\begin{enumerate}
\item $\mu_2(A_0,p_0)>\mu_3(A_0,p_0)$ and $p_0$ is a non-degenerate critical point of
$G_1^-(p)$,
or
\item $\mu_1(A_0,p_0)>\mu_2(A_0,p_0)>\mu_3(A_0,p_0)$ and $p_0$ is a non-degenerate critical point of
$G_2^-(p)$,
or
\item $\mu_1(A_0,p_0)>\mu_2(A_0,p_0)$, $\sqrt{\mu_1(A_0,p_0)}<\sqrt{\mu_2(A_0,p_0)}+\sqrt{\mu_3(A_0,p_0)}$ and
$p_0$ is a non-degenerate critical point of
$G_3^-(p);$
\end{enumerate}
\item $\det M(A_0,p_0)=0$ and
\begin{enumerate}
\item $\mu_2(a_0,p_0)>0$ and $p_0$ is a non-degenerate critical point of
$G_1^0(p)$,
or
\item $\mu_1(A_0,p_0)>\mu_2(A_0,p_0)>0$ and $p_0$ is a non-degenerate critical point of
$G_2^0(p).$
\end{enumerate}
\end{enumerate}
Then, there exists $\epsilon_0>0$ and  a family of connections
$\{{A}_\epsilon\}$ indexed by $\epsilon\in (0, \epsilon_0]$ with the
following properties: ${A}_{\epsilon}$ is a solution to
$\bigl(\mathcal{D}_{\epsilon}\bigr)$ in $\mathcal{A}_{+1}(A_0)$;
$\epsilon^2|{F_{{A}_{\epsilon}}}^{\epsilon}|^2\,dx\to
8\pi^2\delta_{p_0}$ as $\epsilon\to 0$ in the sense of measures
(i.e. $\epsilon {F_{{A}_{\epsilon}}}^{\epsilon}$ concentrates at
$p_0$ as $\epsilon\to 0$).
\end{theorem}

\noindent \textit{Proof:}
Since the different cases can all be proved by very similar arguments, we
only show the proof for the case (1)-(b).

\noindent\emph{Case (1)-(b):} Let $p_0\in
B^4$ satisfy the hypotheses in (1)-(b).  By Lemma \ref{L4.2},
$\langle\Je'(\q),\q_1\rangle=\langle
J'_{\epsilon}(\q),\q_1\rangle+r_{4,1}(\q)=a_{11}(\q)\frac{\partial
J_{\epsilon}}{\partial p_1}(\q)+r_{4,1}(\q).$
Since $a_{11}(\q)\simeq\epsilon^{3/2}$ (cf. Lemma 3.2 in \cite{IM3}),
$|r_{4,1}(\q)|\lesssim\epsilon^4$, and
$J'_{\epsilon}(\q)=\mathcal{F}'_{\epsilon}(\q)+r_2(\q)$ (cf.
Proposition 3.1 in \cite{IM1}), the condition $\langle\Je'(\q),\q_1\rangle=0$ is
equivalent to
\begin{equation}
\frac{\partial\mathcal{F}_{\epsilon}}{\partial
p_1}(\q)+r_{7,1}(\q)=0, \quad \hbox{ with } |r_{7,1}(\q)|\lesssim\epsilon^{5/2}\,.
\end{equation}
More in general, by the same arguments,
$\langle\Je(\q),\q_i\rangle=0$ for $i=1,...,8$, or equivalently $\Je'(\q)=0$, is equivalent to the system
\begin{align}
\label{4.57-4.59}
&\frac{\partial\mathcal{F}_{\epsilon}}{\partial
p_i}(\q)+r_{7,i}(\q)=0\quad\text{with} \quad
|r_{7,i}(\q)|\lesssim\epsilon^{5/2}\quad
\text{for $1\le i\le 4$},\notag\\
&\frac{\partial\mathcal{F}_{\epsilon}}{\partial\xi_i([g])}(\q)+r_{7,i+4}(\q)=0\quad
\text{with}\quad |r_{7,i+4}(\q)|\lesssim\epsilon^{5/2}\quad
\text{for $5\le i\le
7$},\notag\\
&\frac{\partial\mathcal{F}_{\epsilon}}{\partial\lambda}(\q)+r_{7,8}(\q)=0\quad
\text{with} \quad |r_{7,8}(\q)|\lesssim\epsilon^2.
\end{align}
We shall find a solution $\q=(p,[g],\lambda)$  to \eqref{4.57-4.59}, which
satisfies the following conditions:
\begin{enumerate}[(i)]
\item $p$ is in some small neighborhood of $p_0$.
\item $[g]=\exp\xi[g_2^{+}(p)]$,  $\xi\in\mathfrak{so}(3)$ with small $|\xi|$, where $[g_2^+(p)]\in
SO(3)$ is  a critical point for the function
\begin{equation}
\label{4.60}
SO(3)\ni[g]\mapsto\int_{B^4}\bigl((d\underline{A}_0)^-,g(dh_p)^-g^{-1}\bigr)\,dx
\;,
\end{equation}
 with critical value $\Gamma_2^+(A_0,p):=\sqrt{\mu_1(A_0,p)}-\sqrt{\mu_2(A_0,p)}-\sqrt{\mu_3(A_0,p)}$
 (cf. (1)-(b) in $\S3$).

\noindent (Note that, under our hypotheses,
$\Gamma_2^+(A_0,p)>0$ for all $p$ near $p_0$).
\item $\lambda=\lambda_2^+(p)(1+\eta)$, where $|\eta|$ is small and $\lambda_2^+(p)=\Big(\epsilon\frac{\Gamma_2^+(A_0,p)}{F(p)}\Big)^{1/2}$.

\noindent (Note that, by (ii) and Lemma \ref{L4.3} (1), $\lambda_2^+(p)$ is
well-defined if $p$ is close to $p_0$).
\end{enumerate}

We first rewrite the first equation in \eqref{4.57-4.59}. Directly from the definitions above of
$[g_2^+(p)]$ and $\lambda_2^+(p)$, it follows that
\begin{equation}
\label{4.61}
\mathcal{F}_{\epsilon}(p,[g_2^+(p)],\lambda_2^+(p))=-2\epsilon^2G_2^+(p).
\end{equation}
 For $p$ near $p_0$, in the case under consideration one has  $\det M(A_0,p)>0$, $\Gamma^+_2 (A_0, p)>0$ and
 $\mu_1(A_0,p)>\mu_2(A_0,p)>0$. Thus,  the critical point $[g_2^+(p)]$
 is a non-degenerate critical point for the function \eqref{4.60} and, by the implicit function
 theorem, $p\mapsto[g_2^+(p)]$ is differentiable near $p_0$ (cf. $\S3$ (1)-(b)).
From \eqref{4.61}, it follows
\begin{equation}
\label{4.62}
\frac{\partial\mathcal{F}_{\epsilon}}{\partial
p_i}(p,[g_2^+(p)],\lambda_2^+(p))=-2\epsilon^2\frac{\partial
G_2^+}{\partial p_i}(p).
\end{equation}
For $\q=\q(p,\xi,\eta):=(p,\exp\xi[g_2^+(p)],\lambda_2^+(p)(1+\eta))$,
the Taylor's formula yields
\begin{equation}
\label{4.63}
\frac{\partial\mathcal{F}_{\epsilon}}{\partial p_i}(\q(p,\xi,\eta))=
\frac{\partial\mathcal{F}_{\epsilon}}{\partial
p_i}(p,[g_2^+(p)],\lambda_2^+(p))+r_{8,i}(p,\xi,\eta)\quad\text{with $|r_{8,i}(p,\xi,\eta)|\lesssim\epsilon^2(|\xi|+|\eta|)$}.
\end{equation}

\noindent By \eqref{4.62}, \eqref{4.63}, one can rewrite the first equation in \eqref{4.57-4.59} for
$\q=\q(p,\xi,\eta)$ as
\begin{equation}
\label{4.64}
\frac{\partial G_2^+}{\partial
p_i}(p)=r_{9,i}(p,\xi,\eta;\epsilon)\quad\text{with $|r_{9,i}(p,\xi,\eta;\epsilon)|\lesssim |\xi|+|\eta|+\epsilon^{1/2}$}.
\end{equation}
Also, by Taylor's formula,
\begin{align}
\label{4.65}
\frac{\partial G_2^+}{\partial p_i}(p)&=\frac{\partial
G_2^+}{\partial p_i}(p_0)+
\sum_{j=1}^4\frac{\partial^2 G_2^+}{\partial p_i\partial p_j}(p_0)(p_j-(p_0)_j)+r_{10,i}(p)\notag\\
&=\sum_{j=1}^4\frac{\partial^2 G_2^+}{\partial p_i\partial
p_j}(p_0)(p_j-(p_0)_j)+r_{10,i}(p)\quad\text{with $|r_{10,i}(p)|\lesssim|p-p_0|^2$}.
\end{align}
By hypothesis, $p_0$ is a non-degenerate critical point for $G_2^+$,
thus, by \eqref{4.64}, \eqref{4.65}, we may rewrite the first equation of the system \eqref{4.57-4.59} as
\begin{equation}
\label{4.66}
p-p_0=\mathfrak{r}_1(p,\xi,\eta;\epsilon)\quad\text{with $|\mathfrak{r}_1(p,\xi,\eta;\epsilon)|\lesssim|p-p_0|^2+|\xi|+|\eta|+\epsilon^{1/2}$}.
\end{equation}

We next rewrite the second equation in \eqref{4.57-4.59}.

\noindent We denote by $G(p,[g])$ the function \eqref{4.60}. Since
$\frac{\partial\mathcal{F}_{\epsilon}}{\partial\xi_i([g])}(p,[g],\lambda)=-4\epsilon\lambda^2\frac{\partial
G}{\partial\xi_i([g])}(p,[g])$, we have
\begin{align}
\label{4.67}
&\frac{\partial\mathcal{F}_{\epsilon}}{\partial\xi_i([g])}(\q(p,\xi,\eta))=
-4\epsilon\lambda_2^+(p)^2(1+\eta)^2\frac{\partial G}{\partial\xi_i([g])}(p,\exp\xi[g_2^+(p)]) \notag\\
=&~-4\epsilon\lambda_2^+(p)^2(1+\eta)^2\biggl(\frac{\partial
G}{\partial\xi_i([g])}(p,[g_2^+(p)])+ \frac{\partial^2
G}{\partial\xi_i([g])\partial\xi_j([g])}(p,[g_2^+(p)])\xi_j
+r_{10,i+4}(p,\xi)\biggr)\notag\\
=&~-4\epsilon\lambda_2^+(p)^2(1+\eta)^2\biggl(\frac{\partial^2G}{\partial\xi_i([g])\partial\xi_j([g])}(p,[g_2^+(p)])\xi_j+
r_{10,i}(p,\xi)\biggr)\quad\text{with $|r_{10,i}(p,\xi)|\lesssim
|\xi|^2$}.
\end{align}
We now assume $|\eta|<1/2$. Then, by the non-degeneracy of
the critical point $[g_2^+(p)]$ and $\lambda_2^+(p)\simeq\epsilon^{1/2}$ for $p$ near $p_0$, the second equation in \eqref{4.57-4.59} can be  rewritten as
\begin{align}
\label{4.68}
\xi=\mathfrak{r}_2(p,\xi,\eta;\epsilon)\;,\; \text{
with}\quad|\mathfrak{r}_2(p,\xi,\eta;\epsilon)|\lesssim|\xi|^2+\epsilon^{1/2}.
\end{align}

Finally, we rewrite the third equation in \eqref{4.57-4.59}.

\noindent By
$\frac{\partial\mathcal{F}_{\epsilon}}{\partial\lambda}(p,[g],\lambda)=4\lambda^3F(p)-4\epsilon\lambda
G(p,[g])$, one obtains
\begin{align*}
&\frac{\partial\mathcal{F}_{\epsilon}}{\partial\lambda}(\q(p,\xi,\eta))=
4\lambda_2^+(p)(1+\eta)^3F(p)-4\epsilon\lambda_2^+(p)(1+\eta)G(p,\exp\xi[g_2^+(p)])\notag\\
=&~4\lambda_2^+(p)^3F(p)-4\epsilon\lambda_2^+(p)G(p,[g_2^+(p)])+12\lambda_2^+(p)^3\eta
F(p)
-4\epsilon\lambda_2^+(p)\eta G(p,[g_2^+(p)])+r_{10,8}(p,\xi,\eta)\notag\\
=&~(12\lambda_2^+(p)^3F(p)-4\epsilon\lambda_2^+(p)G(p,[g_2^+(p)]))\eta+r_{10,8}(p,\xi,\eta)\notag\\
=&~8\lambda_2^+(p)\epsilon\Gamma_2^+(A_0,p)\eta+r_{10,8}(p,\xi,\eta)\quad\text{with
$|r_{10,8}(p,\xi,\eta)|\lesssim\epsilon^{3/2}(|\xi|^2+|\eta|^2)$}.
\end{align*}
From this, it follows that the last equation in \eqref{4.57-4.59} is equivalent to
\begin{align}
\label{4.70}
\eta=\mathfrak{r}_{3}(p,\xi,\eta;\epsilon)\quad\text{with $|\mathfrak{r}_3(p,\xi,\eta)|\lesssim|\xi|^2+|\eta|^2+\epsilon^{1/2}$}.
\end{align}
We now solve the system given by the equations \eqref{4.66}, \eqref{4.68}, \eqref{4.70} (equivalent to the system \eqref{4.57-4.59}), for small positive $\epsilon$ and
within the set of parameters
$$
M_{\epsilon}:=\{(p,\xi,\eta)\in
B^4\times\mathfrak{so}(3)\times(-1/2,1/2):|p-p_0|
\le\epsilon^{1/8},~|\xi|\le\epsilon^{1/4},~|\eta|\le\epsilon^{1/4}\}.
$$
To this purpose, for $(p,\xi,\eta)\in M_{\epsilon}$, $0\le t\le 1$,
we define the vector function
\begin{equation}
\notag
H_t(p,\xi,\eta;\epsilon):=\left(\begin{array}{c}
p-p_0\\
\xi\\
\eta\\
\end{array}\right)
-t\left(\begin{array}{c}\mathfrak{r}_1(p,\xi,\eta;\epsilon)\\
\mathfrak{r}_2(p,\xi,\eta;\epsilon)\\
\mathfrak{r}_3(p,\xi,\eta;\epsilon)\\
\end{array}\right).
\end{equation}
For $p\in B^4$ with $|p-p_0|=\epsilon^{1/8}$, $\xi\in\mathfrak{so}(3)$ with $|\xi|\le\epsilon^{1/4}$ and $\eta\in\R$ with $|\eta|\le\epsilon^{1/4}$, from \eqref{4.66} it
follows
\begin{equation}
\label{4.72}
|p-p_0|=\epsilon^{1/8}>|\mathfrak{r}_1(p,\xi,\eta;\epsilon)|\,.
\end{equation}
Similarly, for $\xi\in\mathfrak{so}(3)$ with
$|\xi|=\epsilon^{1/4}$, $p\in B^4$ with $|p-p_0|\le\epsilon^{1/8}$ and $\eta\in\R$ with $|\eta|\le\epsilon^{1/4}$, from \eqref{4.68} it follows
\begin{equation}
\label{4.73}
|\xi|=\epsilon^{1/4}>|\mathfrak{r}_2(p,\xi,\eta;\epsilon)|\,
\end{equation}
while,
for $\eta\in\R$ with $\epsilon^{1/4}=\vert\eta\vert$, $p\in B^4$ with $|p-p_0|\le\epsilon^{1/8}$ and $\xi\in\mathfrak{so}(3)$ with $|\xi|\le\epsilon^{1/4}$,   the inequality \eqref{4.70} yields
\begin{equation}
\label{4.74}
|\eta|=\epsilon^{1/4}>|\mathfrak{r}_3(p,\xi,\eta;\epsilon)|\,.
\end{equation}
Thus, from \eqref{4.72}--\eqref{4.74}, it follows that
$H_t(p,\xi,\eta;\epsilon)\ne 0$ for $(p,\xi,\eta)\in\partial
M_{\epsilon}$ and $0\le t\le 1$. Assume by contradiction that
$H_1(p,\xi,\eta;\epsilon)\ne 0$ for all $(p,\xi,\eta)\in
M_{\epsilon}.$ Under this hypothesis, the function
$\frac{H_1}{|H_1|}:M_{\epsilon}\cong B^8\to S^7$ would be
well-defined and continuous, hence homotopically trivial when
restricted to $\partial M_{\epsilon}$. On the other hand, the
well-defined maps $\frac{H_t}{|H_t|}:\partial M_{\epsilon}\cong
S^7\to S^7$, for $0\le t\le 1$, give a homotopy between
$\frac{H_0}{|H_0|}$ and $\frac{H_1}{|H_1|}$. Since
$\frac{H_0}{|H_0|}:\partial M_{\epsilon}\cong S^7\to S^7$ is
obviously homotopically non-trivial (in fact, its topological degree
is $1$), we have a contradiction. Thus, there exists
$(p,\xi,\eta)\in M_{\epsilon}$ such that
$H_1(p,\xi,\eta;\epsilon)=0$, i.e., a solution of the system given by the equations \eqref{4.66}, \eqref{4.68}, \eqref{4.70}.
Since this solution satisfies
$|p-p_0|\le\epsilon^{1/8}\to 0$ as $\epsilon\to 0$, it concentrates
at $p_0$ as $\epsilon\to 0$. This completes the proof of Theorem 2. \hfill$\Box$

\section{Proof of Theorem 3}
The third and final existence theorem is the following.
\begin{theorem}
\label{3}
Assume that there exists $p_0\in B^4$ such that one of the following holds:

\begin{enumerate}[(1)]
\item   $\det M(A_0,p_0)>0$,
$\mu_1(A_0,p_0)>\mu_2(A_0,p_0)>\mu_3(A_0,p_0)$ and
$\sqrt{\mu_1(A_0,p_0)}>\sqrt{\mu_2(A_0,p_0)}+\sqrt{\mu_3(A_0,p_0)}$;
\item $\det M(A_0,p_0)<0$ and $\mu_1(A_0,p_0)>\mu_2(A_0,p_0)>\mu_3(A_0,p_0)$;
\item $\det M(A_0,p_0)=0$ and $\mu_1(A_0,p_0)>\mu_2(A_0,p_0)>0$;
\end{enumerate}
then, for all sufficiently small $\epsilon>0$, there exist at least
two distinct solutions to $(\mathcal{D}_{\epsilon})$ in
$\mathcal{A}_{+1}(A_0).$
Furthermore, the following
alternative holds: there
exists at least one non-minimizing solution, or there exist
infinitely many minimizing solutions. In the hypotheses (2), if in addition $\sqrt{\mu_1(A_0,p_0)}<\sqrt{\mu_2(A_0,p_0)}+\sqrt{\mu_3(A_0,p_0)},$
then  there exist at least three distinct solutions, of which at least two
non-minimizing, or there exist infinitely many minimizing solutions
to $(\mathcal{D}_{\epsilon})$ in $\mathcal{A}_{+1}(A_0)$.
\end{theorem}
In order to prove Theorem \ref{3} we need to prove Lemma \ref{L4.10}, which enables us to apply the standard critical point theory (by showing that a subset of the parameter space, namely  the set $\Je^{C_{\epsilon}-C_0\epsilon^2}(d_0,D_1,D_2)$ defined below, is invariant under the negative gradient flow of $\Je$),  Lemma \ref{L4.11} and Corollary \ref{C4.1}, where the topological properties of
$\Je^{C_{\epsilon}-C_0\epsilon^2}(d_0,D_1,D_2)$ are studied.

To this purpose, let $C_{\epsilon}:=8\pi^2+\int_{B^4}|{F_{\underline{A}_{\epsilon}}}^{\epsilon}|^2\,dx$ be the constant in \eqref{4.50}, and let us define
\begin{equation}
\label{4.75}
\Je^{C_{\epsilon}-C_0\epsilon^2}(d_0,D_1,D_2)=\{\q=(p,[g],\lambda)\in
B_{1-d_0}\times SO(3)\times(\sqrt{D_1\epsilon},\sqrt{D_2\epsilon}):
\Je(\q)\le C_{\epsilon}-C_0\epsilon^2\}\,.
\end{equation}
\begin{lemma}
\label{L4.10} Let $C_0$ be any given positive constant.
There exist $0<d_0<1$, $0<D_1<D_2$ and $\epsilon_0>0$ such that the following properties hold for all $\epsilon$ satisfying $0<\epsilon<\epsilon_0$:
\begin{equation}
\label{4.76}
\Big\langle\frac{\partial\Je}{\partial
p}(\q),\frac{p}{|p|}\Big\rangle>0, \hbox{ for }\q\in\Je^{C_{\epsilon}-C_0\epsilon^2}(d_0,D_1,D_2),\;
p\in\partial B_{1-d_0}\,,
\end{equation}
\begin{equation}
\label{4.77}
\frac{\partial\Je}{\partial\lambda}(\q)<0, \hbox{ for }\q\in\Je^{C_{\epsilon}-C_0\epsilon^2}(d_0,D_1,D_2),\;
\lambda=\sqrt{D_1\epsilon}\,,
\end{equation}
\begin{equation}
\label{4.78}
\frac{\partial\Je}{\partial\lambda}(\q)>0, \hbox{ for }\q\in\Je^{C_{\epsilon}-C_0\epsilon^2}(d_0,D_1,D_2),\;
\lambda=\sqrt{D_2\epsilon}\,.
\end{equation}
\end{lemma}
\textit{Proof:} By Lemma \ref{L4.1} and Proposition 3.1 in \cite{IM1},
\begin{equation}
\notag
C_{\epsilon}+\mathcal{F}_{\epsilon}(\q)+r_1(\q)+r_3(\q)\le
C_{\epsilon}-C_0\epsilon^2,
\end{equation}
for $\q\in\PP(d_0,\lambda_0;D_1,D_2;\epsilon)$ with $\Je(\q)\le
C_{\epsilon}-C_0\epsilon^2$.

Since $|r_1(\q)|\lesssim\epsilon^3$ and
$|r_3(\q)|\lesssim\epsilon^3$, there exists
$\epsilon(d_0,D_1,D_2)>0$ (depending on $d_0$, $D_1$ and $D_2$) such
that
\begin{equation}
\label{4.80}
\mathcal{F}_{\epsilon}(\q)\le-\frac{C_0}{2}\epsilon^2
\end{equation}
for $0<\epsilon<\epsilon(d_0,D_1,D_2)$ and $\q\in\Je^{C_{\epsilon}-C_0\epsilon^2}(d_0,D_1,D_2)$.

\noindent From (1) of Lemma \ref{L4.3},
$-4\epsilon\lambda^2G(p,[g])<\mathcal{F}_{\epsilon}(\q)$ and from \eqref{4.80}
it follows that
\begin{equation}
\label{4.81}
G(p,[g])\ge\frac{C_0}{8}\epsilon\lambda^{-2}
\end{equation}
for
$\q=(p,[g],\lambda)\in\Je^{C_{\epsilon}-C_0\epsilon^2}(d_0,D_1,D_2)$,
if $0<\epsilon<\epsilon(d_0,D_1,D_2)$.

On the other hand, by (3) of Lemma \ref{L4.3}  and Lemma \ref{L4.4}, one has
\begin{align}
\notag
\frac{\partial\mathcal{F}_{\epsilon}}{\partial
p_i}(\q)&=2\lambda^4\frac{\partial F}{\partial p_i}(p)-
4\epsilon\lambda^2\frac{\partial G}{\partial p_i}(p,[g]) \ge
C_2D_1^2\epsilon^2\frac{1}{d_0^5}\frac{p_i}{|p|}-4\epsilon^2D_2C_3
\end{align}
for $\q\in\Je^{C_{\epsilon}-C_0\epsilon^2}(d_0,D_1,D_2)$ with
$p\in\partial B_{1-d_0}$ for some absolute constants $C_2, C_3 >0.$

\noindent This inequality implies
\begin{equation}
\label{4.83}
\Big\langle\frac{\partial\mathcal{F}_{\epsilon}}{\partial
p},\frac{p}{|p|}\Big\rangle\ge
C_2D_1^2\epsilon^2\frac{1}{d_0^5}-
4\epsilon^2D_2C_3.
\end{equation}
From Proposition 3.1 in \cite{IM1}, Lemma \ref{L4.2} and Lemma 3.2 in \cite{IM3} (cf. also the proof of Theorem
\ref{2} in the previous section), one obtains
\begin{equation}
\label{4.84}
\frac{\partial\Je}{\partial
p_i}(\q)=\frac{\partial\mathcal{F}_{\epsilon}}{\partial
p_i}(\q)+r_{11,i}(\q),\quad |r_{11,i}(\q)|
\lesssim\epsilon^{5/2}~\quad(1\le i\le 4)
\end{equation}
and, combining \eqref{4.83}, \eqref{4.84},
\begin{equation}
\label{4.85}
\Big\langle\frac{\partial\Je}{\partial
p}(\q),\frac{p}{|p|}\Big\rangle\ge
C_2D_1^2\epsilon^2\frac{1}{d_0^5}- 4\epsilon^2D_2C_3+r_{12}(\q)
\end{equation}
for $\q\in\Je^{C_{\epsilon}-C_0\epsilon^2}(d_0,D_1,D_2)$ with
$p\in\partial B_{1-d_0}$, where
$|r_{12}(\q)|\lesssim\epsilon^{5/2}$.

We next estimate
$\frac{\partial\mathcal{F}_{\epsilon}}{\partial\lambda}(\q)=8\lambda^3F(p)-8\epsilon\lambda
G(p,[g])$.

\noindent From \eqref{4.81}, one obtains
\begin{align}
\label{4.86}
\frac{\partial\mathcal{F}_{\epsilon}}{\partial\lambda}(\q)&=8\sqrt{D_1\epsilon}(D_1\epsilon
F(p)-\epsilon G(p,[g])) \le 8\sqrt{D_1\epsilon}\Big(D_1\epsilon
F(p)-\frac{C_0}{8D_1}\epsilon\Big)\,,
\end{align}
for $\q\in\Je^{C_{\epsilon}-C_0\epsilon^2}(d_0,D_1,D_2)$,
$\lambda=\sqrt{D_1\epsilon}$.

\noindent Similarly, by Lemma \ref{L4.4},
\begin{equation}
\label{4.87}
\frac{\partial\mathcal{F}_{\epsilon}}{\partial\lambda}(\q)\ge8\sqrt{D_2\epsilon}(D_2\epsilon
F(p)-C_4\epsilon)
\end{equation}
for $\q\in\Je^{C_{\epsilon}-C_0\epsilon^2}(d_0,D_1,D_2)$,
$\lambda=\sqrt{D_2\epsilon}$, for some absolute constant $C_4>0$.

We now observe that, by (1), (2) of Lemma \ref{L4.3}, there exists an
absolute constant $C_5>0$ such that $F(p)\ge C_5$ for all $p\in
B^4$. To complete the proof, we choose $D_2>0$ such that
$D_2C_5-C_4>\frac{C_4}{2}$ (take for example $D_2=\frac{2C_4}{C_5}$
- notice that this is independent of $d_0$), and $D_1=Dd_0^2$, where
the absolute constant $D>0$ is chosen to satisfy
\begin{equation}
\label{4.88}
Dd_0^2F(p)-\frac{C_0}{8Dd_0^2}<-\frac{C_0}{16Dd_0^2}\quad\text{for all
$p\in B_{1-d_0}$},
\end{equation}
i.e., $D^2d_0^4F(p)<\frac{C_0}{16}$ holds for all $p\in B_{1-d_0}$.
(By (2) of Lemma \ref{L4.3}, one can easily see that such constant $D$
exists).

\noindent Furthermore, we choose $0<d_0<1$ so that
\begin{equation}
\label{4.89}
C_2D^2\frac{1}{d_0}-4D_2C_3\ge 2.
\end{equation}
By \eqref{4.85}, there exists $0<\epsilon_0<\epsilon(d_0,D_1,D_2)$ such
that for $0<\epsilon<\epsilon_0$ and
$\q\in\Je^{C_{\epsilon}-C_0\epsilon^2}(d_0,D_1,D_2)$,
$p\in\partial B_{1-d_0}$, one has
\begin{equation}
\label{4.90}
\Big\langle\frac{\partial\Je}{\partial
p},\frac{p}{|p|}\Big\rangle\ge 2\epsilon^2+r_{12}(\q)>0.
\end{equation}
From \eqref{4.86}, \eqref{4.88}, Proposition 3.1 in \cite{IM1}, and Lemma \ref{L4.2}  (cf. also
the proof of Theorem 2), we also have
\begin{align}
\frac{\partial\Je}{\partial\lambda}(\q)=\frac{\partial\mathcal{F}_{\epsilon}}{\partial\lambda}(\q)+r_{13}(\q)
\le-\frac{C_0}{2\sqrt D_1 d_0}\epsilon^{3/2}+r_{13}(\q)<0, \qquad (\hbox{with  }|r_{13}(\q)|
\lesssim\epsilon^2)\,,
\end{align}
for $0<\epsilon<\epsilon_0$ and
$\q\in\Je^{C_{\epsilon}-C_0\epsilon^2}(d_0,D_1,D_2)$, with $\lambda=\sqrt{D_1\epsilon}$, if
$\epsilon_0>0$ is chosen suitably
small. Similarly, by \eqref{4.87},
\begin{align}
\label{4.92}
\frac{\partial\Je}{\partial\lambda}(\q)=\frac{\partial\mathcal{F}_{\epsilon}}{\partial\lambda}(\q)+r_{13}(\q)
\ge 4 \sqrt D_2 C_4\epsilon^{3/2}+r_{13}(\q)>0,  \hbox{ with }\lambda=\sqrt{D_2\epsilon}
\end{align}
for $0<\epsilon<\epsilon_0$ and
$\q\in\Je^{C_{\epsilon}-C_0\epsilon^2}(d_0,D_1,D_2)$,
 if $\epsilon_0>0$ is chosen suitably
small.
\noindent The proof now follows from \eqref{4.90} -- \eqref{4.92}. \hfill$\Box$

\medskip

\noindent As a consequence of Lemma \ref{L4.10}, the set
$\Je^{C_{\epsilon}-C_0\epsilon^2}(d_0,D_1,D_2)$ is invariant under
the negative gradient flow of $\Je$, i.e., the solution of the
differential equation $\frac{d\q(t)}{dt}=-\Je'(\q(t))$ stays in
$\Je^{C_{\epsilon}-C_0\epsilon^2}(d_0,D_1,D_2)$ for all $t\ge 0$.
Therefore, one can apply the standard critical point theory
(see~\cite{Chang},~\cite{Schwartz},~\cite{Struwe}) to the function
$\Je$ restricted on $\Je^{C_{\epsilon}-C_0\epsilon^2}(d_0,D_1,D_2)$.

The number of critical points of $\Je$ on
$\Je^{C_{\epsilon}-C_0\epsilon^2}(d_0,D_1,D_2)$ depends on the
topological complexity of
$\Je^{C_{\epsilon}-C_0\epsilon^2}(d_0,D_1,D_2)$. To study the topological
properties of this set, we consider the function
$G(p,[g])=\int_{B^4}\bigl(F_{\underline{A}_0}^-,g(dh_p)^-g^{-1}\bigr),$ and
for $p_0\in B^4$,
$\eta>0$, we take the following subset of $SO(3)$:
\begin{equation}
\label{4.93}
\mathcal{S}(p_0,\eta)=\{[g]\in SO(3):-G(p_0,[g])\le-\eta\}\;.
\end{equation}

We recall the following topological notion
(see~\cite{Chang},~\cite{Schwartz},~\cite{Struwe} for more details):
let $X$ be a topological space. The Ljusternik-Schnirelman category (LS-category in short) of a closed subset
$A\subset X$ with respect to $X$, denoted by $\cat(A,X)$, is defined
as the least integer $k$ such that $A\subset A_1\cup\cdots\cup A_k$,
where the $A_i$'s (for $i=1,2,\ldots,k$) are closed and contractible
in $X$.

\noindent We have the following:
\begin{lemma} \label{L4.11} Assume that there exists $p_0\in B^4$ such that
$\mu_1(A_0,p_0)>\mu_2(A_0,p_0)>\mu_3(A_0,p_0)$. Then there exists $\eta>0$ such that the following holds:
\begin{enumerate}[(1)]
\item
If $p_0$ satisfies $\det M(A_0,p_0)\ge 0$ and $\sqrt{\mu_1(A_0,p_0)}>\sqrt{\mu_2(A_0,p_0)}+\sqrt{\mu_3(A_0,p_0)}$, then
$\cat(\mathcal{S}(p_0,\eta),SO(3))\ge 2$.
\item If $p_0$ satisfies
$\det M(A_0,p_0)<0$, then $\cat(\mathcal{S}(p_0,\eta),SO(3))\ge 2$.
If in addition $p_0$ satisfies
$\sqrt{\mu_1(A_0,p_0)}<\sqrt{\mu_2(A_0,p_0)}+\sqrt{\mu_3(A_0,p_0)}$,
then $\cat(\mathcal{S}(p_0,\eta),SO(3))\geq 3$.
\end{enumerate}
\end{lemma}
\textit{Proof:} The proof follows from Morse theory
(cf.~\cite{Milnor}). If
$\mu_1(A_0,p_0)>\mu_2(A_0,p_0)>\mu_3(A_0,p_0)$, the function
$G(p_0,\cdot)$ on $SO(3)$ is a Morse function with four critical
points, with Morse indexes equal to $3$, $2$, $1$, and $0$ (cf.
$\S3$). By Morse theory, this yields a cell decomposition of
$SO(3)$ as follows
\begin{equation}
\label{4.94}
SO(3)\cong e^0\cup e^1\cup e^2\cup e^3,
\end{equation}
where $e^i$ ($i=0,1,2,3$) is a cell of dimension $i$. We examine the
two cases separately.

\textit{(1) $\det M(A_0,p_0)\ge 0$,
$\sqrt{\mu_1(A_0,p_0)}>\sqrt{\mu_2(A_0,p_0)}+\sqrt{\mu_3(A_0,p_0)}$.}

\noindent In this case, two of the critical points of the function
$G(p_0,\cdot)$ on $SO(3)$ assume positive critical values, namely
$\sqrt{\mu_1(A_0,p_0)}+\sqrt{\mu_2(A_0,p_0)}+\sqrt{\mu_3(A_0,p_0)}$
and
$\sqrt{\mu_1(A_0,p_0)}-\sqrt{\mu_2(A_0,p_0)}-\sqrt{\mu_3(A_0,p_0)}$,
with Morse indexes 3 and 2, respectively (cf. (1) of $\S3$). Set
$\eta=\frac{1}{2}(\sqrt{\mu_1(A_0,p_0)}-\sqrt{\mu_2(A_0,p_0)}-\sqrt{\mu_3(A_0,p_0)})$.
By Morse theory, $\mathcal{S}(p_0,\eta)$ is homotopically equivalent
to $e^0\cup e^1$, we write $\mathcal{S}(p_0,\eta)\cong e^0\cup e^1$.
(This is meant with respect to the cell decomposition `dual' to
\eqref{4.94}). Since $H_1(SO(3);\Z_2)\simeq\Z_2$ is generated by the cell
$e^1$, $\mathcal{S}(p_0,\eta)\subset SO(3)$ is not contractible.
This yields $\cat(\mathcal{S}(p_0,\eta),SO(3))\ge 2$.

\textit{(2) $\det M(A_0,p_0)<0$.}

\noindent In this case, at least two of the critical points of the
function $G(p_0,\cdot)$ on $SO(3)$ assume positive critical values,
namely
$\sqrt{\mu_1(A_0,p_0)}+\sqrt{\mu_2(A_0,p_0)}-\sqrt{\mu_3(A_0,p_0)}$
and
$\sqrt{\mu_1(A_0,p_0)}-\sqrt{\mu_2(A_0,p_0)}+\sqrt{\mu_3(A_0,p_0)}$.
Setting $\eta=\sqrt{\mu_3(A_0,p_0)}$, one can show by the same
argument as in the previous case, that the set
$\mathcal{S}(p_0,\eta)$ is not contractible in $SO(3)$ and
$\cat(\mathcal{S}(p_0,\eta),SO(3))\ge 2$. In the additional
hypothesis
$\sqrt{\mu_1(A_0,p_0)}<\sqrt{\mu_2(A_0,p_0)}+\sqrt{\mu_3(A_0,p_0)}$,
three of the critical points of $G(p_0,\cdot)$ on $SO(3)$ assume
positive critical values, namely
$\sqrt{\mu_1(A_0,p_0)}+\sqrt{\mu_2(A_0,p_0)}-\sqrt{\mu_3(A_0,p_0)}$,
$\sqrt{\mu_1(A_0,p_0)}-\sqrt{\mu_2(A_0,p_0)}+\sqrt{\mu_3(A_0,p_0)}$
and
$-\sqrt{\mu_1(A_0,p_0)}+\sqrt{\mu_2(A_0,p_0)}+\sqrt{\mu_3(A_0,p_0)}$.
Set
$\eta=\frac{1}{2}(-\sqrt{\mu_1(A_0,p_0)}+\sqrt{\mu_2(A_0,p_0)}+\sqrt{\mu_3(A_0,p_0)})$.
By Morse theory, $\mathcal{S}(p_0,\eta)\simeq e^0\cup e^1\cup e^2$ (homotopically
equivalent with respect to the cell decomposition dual to \eqref{4.94}). By
contradiction, assume $\cat(\mathcal{S}(p_0,\eta),SO(3))\le 2$. By
Morse theory, $SO(3)$ is obtained by attaching the $3$-cell $e^3$, in
correspondence with the maximum of $-G(p_0,\cdot)$ and which is contractible in $SO(3)$, to $\mathcal{S}(p_0,\eta)$. Thus, one would
also have $\cat (SO(3),SO(3))\le 3$. On the other hand,
$\cat(SO(3),SO(3))\ge\text{cuplength}\,H^{\ast}(SO(3);\Z_2)+1=4$,
where $\text{cuplength}\,H^{\ast}(SO(3);\Z_2)$ is the cuplength of
the cohomology ring $H^{\ast}(SO(3);\Z_2)\cong\Z_2[a]/(a^4)$, where
$\deg a=1$ (cf.~\cite{Chang},~\cite{Schwartz}). This yields a
contradiction. Thus, $\cat(\mathcal{S}(p_0,\eta),SO(3))\geq 3$.

This completes the proof. \hfill$\Box$

\medskip

Associated with the set $\mathcal{S}(p_0,\eta)$, we define the set
\begin{equation}
\label{4.95}
\tilde{\mathcal{S}}(p_0,\eta):=\{p_0\}\times\mathcal{S}(p_0,\eta)\times\{\lambda_0\}\subset
B_{1-d_0}\times
SO(3)\times(\sqrt{D_1\epsilon},\sqrt{D_2\epsilon})\;,
\end{equation}
where $\lambda_0:=\Big(\frac{\eta\epsilon}{F(p_0)}\Big)^{1/2}$.

The following is a corollary of Lemma \ref{L4.11}.
\newtheorem{corollary}{Corollary}[section]
\begin{corollary} \label{C4.1} In the hypotheses of Lemma \ref{L4.11}, the following
assertions hold for the cases (1) and (2), respectively:
\begin{enumerate}[(1)]
\item $\cat(\tilde{\mathcal{S}}(p_0,\eta),B_{1-d_0}\times SO(3)\times(\sqrt{D_1\epsilon},\sqrt{D_2\epsilon}))\ge 2$.
\item $\cat(\tilde{\mathcal{S}}(p_0,\eta),B_{1-d_0}\times SO(3)\times(\sqrt{D_1\epsilon},\sqrt{D_2\epsilon}))\ge 2$.
Moreover, if $p_0$ satisfies the additional hypothesis
$\sqrt{\mu_1(A_0,p_0)}<\sqrt{\mu_2(A_0,p_0)}+\sqrt{\mu_3(A_0,p_0)}$,
then $$\cat(\tilde{\mathcal{S}}(p_0,\eta),B_{1-d_0}\times
SO(3)\times(\sqrt{D_1\epsilon},\sqrt{D_2\epsilon}))\ge 3.$$
\end{enumerate}
\end{corollary}
\textit{Proof:} (1) Assume, by contradiction, that
$\tilde{\mathcal{S}}(p_0,\eta)$ be contractible in $B_{1-d_0}\times
SO(3)\times(\sqrt{D_1\epsilon},\sqrt{D_2\epsilon})$. Then there
exists $h:[0,1]\times\tilde{\mathcal{S}}(p_0,\eta)\to
B_{1-d_0}\times SO(3)\times (\sqrt{D_1\epsilon},\sqrt{D_2\epsilon})$
such that $h(0,\cdot)=\tilde\iota$, where $\tilde \iota$ is the
inclusion $\tilde \iota:\tilde{\mathcal{S}}(p_0,\eta)\hookrightarrow
B_{1-d_0}\times SO(3)\times
(\sqrt{D_1\epsilon},\sqrt{D_2\epsilon})$,
 and
$h(1,\cdot):=\text{const}$. We define the map $H$ as the following
composition of maps
$$H:[0,1]\times\mathcal{S}(p_0,\eta)\hookrightarrow[0,1]\times
\tilde{\mathcal{S}}(p_0,\eta)\overset{h}{\longrightarrow}B_{1-d_0}\times
SO(3) \times
(\sqrt{D_1\epsilon},\sqrt{D_2\epsilon})\overset{\text{pr}_2}{\longrightarrow}SO(3),$$
where the first inclusion is given by
$(t,[g])\mapsto(t,(p_0,[g],\lambda_0))$ and $\text{pr}_2$ is the
projection on the second factor. The map $H$ satisfies
$H(0,\cdot)=\iota$, where $\iota:
\mathcal{S}(p_0,\eta)\hookrightarrow SO(3)$, and
$H(1,\cdot):=\text{const}$, yielding a contradiction since
$\mathcal{S}(p_0,\eta)$ is not contractible in $SO(3)$.

(2) The first part of the statement is proved as in (1). We need to
prove the second part. By contradiction, assume that, in the given
additional hypothesis,
$\cat(\tilde{\mathcal{S}}(p_0,\eta),B_{1-d_0}\times SO(3)\times
(\sqrt{D_1\epsilon},\sqrt{D_2\epsilon}))=2$. Then there exist
contractible closed sets $\tilde{A}_1,\tilde{A}_2\subset
B_{1-d_0}\times SO(3)\times (\sqrt{D_1\epsilon},\sqrt{D_2\epsilon})$
such that
$\tilde{\mathcal{S}}(p_0,\eta)\subset\tilde{A}_1\cup\tilde{A}_2$.
Let us take the projection $\text{pr}_{1,3}:B_{1-d_0}\times
SO(3)\times(\sqrt{D_1\epsilon},\sqrt{D_2\epsilon})\to
B_{1-d_0}\times(\sqrt{D_1\epsilon},\sqrt{D_2\epsilon})$,  and define
$A_1:=\text{pr}_{1,3}^{-1}(\{(p_0,\lambda_0)\})\cap\tilde{A}_1$ and
$A_2:=\text{pr}_{1,3}^{-1}(\{(p_0,\lambda_0)\})\cap\tilde{A}_2$.
Under the natural identification $SO(3)\cong\{p_0\}\times
SO(3)\times\{\lambda_0\}$, $A_1$ and $A_2$ are closed subsets of
$SO(3)$ and $\mathcal{S}(p_0,\eta)\subset A_1\cup A_2$. Moreover,
the composition of maps
$$[0,1]\times A_1\hookrightarrow[0,1]\times \tilde{A}_1\overset{h_1}{\longrightarrow} B_{1-d_0}
\times SO(3)\times
(\sqrt{D_1\epsilon},\sqrt{D_2\epsilon})\overset{\text{pr}_2}{\longrightarrow}SO(3),$$
with $h_1:[0,1]\times\tilde{A}_1\to B_{1-d_0}\times SO(3)\times
(\sqrt{D_1\epsilon},\sqrt{D_2\epsilon})$ satisfying
$h_1(0,\cdot)=\tilde\iota_1$ (where
$\tilde\iota_1:\tilde{A}_1\hookrightarrow B_{1-d_0}\times
SO(3)\times (\sqrt{D_1\epsilon},\sqrt{D_2\epsilon})$ is the inclusion),
and $h_1(1,\cdot):=\text{const}$, shows that $A_1$ is contractible
in $SO(3)$. Similarly, $A_2$ is contractible in $SO(3)$. This yields
$\cat(\mathcal{S}(p_0,\eta),SO(3))\le 2$, thus a contradiction. This
completes the proof. \hfill$\Box$

\medskip

We are now ready to prove Theorem 3.

\textit{Proof of Theorem 3:} Let $\eta>0$ be as in Lemma \ref{L4.11}. We
first show that for $C_0=\frac{\eta^2}{F(p_0)}$ and for all small
$\epsilon>0$,
\begin{equation}
\label{4.96}
\tilde{\mathcal{S}}(p_0,\eta)\subset\Je^{C_{\epsilon}-C_0\epsilon^2}(d_0,D_1,D_2).
\end{equation}
In fact, for $(p_0,[g],\lambda_0)\in\tilde{\mathcal{S}}(p_0,\eta)$,
\begin{align}
\notag
\mathcal{F}_{\epsilon}(p_0,[g],\lambda_0)&=2\lambda_0^4F(p_0)-4\epsilon\lambda_0^2G(p_0,[g])
\le2\lambda_0^4F(p_0)-4\epsilon\lambda_0^2\eta\le-2\frac{\eta^2}{F(p_0)}\epsilon^2.
\end{align}
From Proposition 3.1 in \cite{IM1} and Lemma \ref{L4.1} it follows that
\begin{align}
\notag
\Je(p_0,[g],\lambda_0)&=C_{\epsilon}+2\mathcal{F}_{\epsilon}(p_0,[g],\lambda_0)+r_1(p_0,[g],\lambda_0)+
r_3(p_0,[g],\lambda_0)\notag\\
&\le
C_{\epsilon}-2\frac{\eta^2}{F(p_0)}\epsilon^2+r_1(p_0,[g],\lambda_0)+r_3(p_0,[g],\lambda_0)\le
C_{\epsilon}-\frac{\eta^2}{F(p_0)}\epsilon^2
\end{align}
if $\epsilon>0$ is small enough.

\noindent By Corollary \ref{C4.1}, we see that in all the different cases in Theorem 3,
the set
\begin{align*}
\mathcal{K}_2:= \{\text{compact sets }&A\subset B_{1-d_0}\times
SO(3)\times (\sqrt{D_1\epsilon},\sqrt{D_2\epsilon}):\notag\\
&\cat(A,B_{1-d_0}\times
SO(3)\times(\sqrt{D_1\epsilon},\sqrt{D_2\epsilon})\ge 2\}
\end{align*}
is non-empty (in fact,
$\tilde{\mathcal{S}}(p_0,\eta)\in\mathcal{K}_2$).

\noindent Let
\begin{equation}
\notag
c_2:=\inf_{A\in\mathcal{K}_2}\{\max\Je(\q):\q\in A\}.
\end{equation}
Then by \eqref{4.96},
\begin{equation}
\label{4.100}
c_2\le C_{\epsilon}-C_0\epsilon^2.
\end{equation}

Combining \eqref{4.100} with Lemma \ref{L4.10}, standard arguments show that
$c_2$ is a critical value for $\Je$. By contradiction, assume not.
Then, by flowing along the negative gradient of
$\Je$, one sees that there exists a $\delta$,
$0<\delta<C_0\epsilon^2$, and a deformation $\Theta$ of
$\Je^{C_{\epsilon}-C_0\epsilon^2}$ (i.e., $\Theta=\theta(1,\cdot)$,
for some continuous
$\theta:[0,1]\times\Je^{C_{\epsilon}-C_0\epsilon^2}\to\Je^{C_{\epsilon}-C_0\epsilon^2}$
with $\theta(0,\cdot)=\text{Id}$) which satisfies
\begin{equation}
\label{4.101}
\Theta(\Je^{c_2+\delta})\subset\Je^{c_2-\delta}.
\end{equation}
Take $A\in\mathcal{K}_2$ such that $\max_{\q\in
A}\Je(\q)<c_2+\delta$. Since $\Theta(A)\in\mathcal{K}_2$
(cf.~\cite{Struwe}) and $\max_{\q\in\Theta(A)}\Je(\q)\le c_2-\delta$
by \eqref{4.101}, this contradicts the definition of $c_2$. Thus $c_2$ is
a critical value of $\Je$.

\noindent Furthermore, by Lemma \ref{L4.10},
$$c_1:=\min\{\Je(\q):\q\in B_{1-d_0}\times SO(3)\times (\sqrt{D_1\epsilon},\sqrt{D_2\epsilon})\}$$
is also a critical value of $\Je$.

If $c_1<c_2$, there are at least two critical points of $\Je$ on
$B_{1-d_0}\times
SO(3)\times(\sqrt{D_1\epsilon},\sqrt{D_2\epsilon})$,  and at least
one of them is a non-minimal critical point. If instead $c_1=c_2$,
then the LS-category of the critical set
$$\mathcal{C}(c_1):=\{q\in
B_{1-d_0}\times
SO(3)\times(\sqrt{D_1\epsilon},\sqrt{D_2\epsilon}):\Je'(\q)=0,~\Je(\q)=c_1\}$$
is at least $2$ (see~\cite{Struwe}), hence,
 $\#\mathcal{C}(c_1)=+\infty$. If, furthermore,
$c_1=\inf_{A\in\mathcal{A}_{+1}(A_0)}\epsilon^2\YMe(A)$, there are
infinitely many minimizing solutions in $\mathcal{A}_{+1}(A_0)$.
Otherwise, i.e., if
$c_1>\inf_{A\in\mathcal{A}_{+1}(A_0)}\epsilon^2\YMe(A)$, there are
infinitely many non-minimal critical points. This complete the
proof, except for the second statement in (2).

To prove the latter, we consider the set
\begin{align*}
\mathcal{K}_3:= \{\text{compact sets }&A\subset B_{1-d_0}\times
SO(3)\times(\sqrt{D_1\epsilon},\sqrt{D_2\epsilon}):\notag\\
&\cat(A,B_{1-d_0}\times
SO(3)\times(\sqrt{D_1\epsilon},\sqrt{D_2\epsilon})\ge 3\},
\end{align*}
which again is non-empty by Corollary \ref{C4.1}, and define
$$c_3:=\inf_{A\in\mathcal{K}_3}\{\max\Je(\q):\q\in A\}\le C_{\epsilon}-C_0\epsilon^2\;,$$
which is a critical value for $\Je$. We may assume  $c_1<c_2$ (if
$c_1 = c_2$, we have already showed that there are infinitely many
solutions). If also $c_2<c_3$, then there are at least two
non-minimizing solutions. If $c_2=c_3$, again, the LS-category of
$\mathcal{C}(c_2)=\{\q\in B_{1-d_0}\times
SO(3)\times(\sqrt{D_1\epsilon},\sqrt{D_2\epsilon}):\Je'(\q)=0,~\Je(\q)=c_2\}$
is at least 2, hence $\#\mathcal{C}(c_2)=+\infty$. Thus, there are
infinitely many non-minimizing solutions. So, the assertion (2) of
Theorem 3 is completely proved for all the different cases.
\hfill$\Box$

\medskip

The proof above yields the following corollary:
\begin{corollary}
\label{C4.2}
Assume there exists a $p_0\in B^4$ and $\eta>0$ such that
$\mathcal{S}(p_0,\eta)$ is not contractible in $SO(3)$. Then there
exist at least two solutions of $(\mathcal{D}_{\epsilon})$ in
$\mathcal{A}_{+1}(A_0)$. Furthermore, if one assumes that
$\cat(\mathcal{S}(p_0,\eta),SO(3))\ge 3$, then there exist at least
three solutions of $(\mathcal{D}_{\epsilon})$ in
$\mathcal{A}_{+1}(A_0)$.
\end{corollary}

\section{Examples}
In this section we perform two tasks. First we illustrate a method
to construct boundary data that yield any prescribed matrix $M$
(precisely, for any given matrix $M$, a boundary value is
constructed such that $M(A_0, p_0)= M$), thus showing that the
different cases in Theorem 3 can all be achieved (Proposition
\ref{P4.1}) . Second, We show that
the non-degeneracy condition
$\mu_1(A_0,p_0)>\mu_2(A_0,p_0)>\mu_3(A_0,p_0)$ (where the
$\mu_i(A_0,p_0)$'s ($i=1,2,3$) are the eigenvalues of $M(A_0, p_0)^t
M(A_0, p_0)$) can always be obtained by making an an arbitrarily
small perturbation of the boundary value $A_0$ (Proposition
\ref{P4.2}). This can be done while leaving the other conditions in
Theorems 1-3 unaffected.

Let $A_0\in C^{\infty}(T^{\ast}\partial B^4\otimes\so(3))$ be a
given boundary value and let $\underline{A}_0\in
C^{\infty}(T^{\ast}B^4\otimes\so(3))$ be any solution to
$(\mathcal{D}_0)$ with boundary value $A_0$ (recall that
$d\underline{A}_0$ is uniquely determined by the boundary data
$A_0$).

Let $\omega_1^-:=dx^1\wedge dx^2-dx^3\wedge dx^4$,
$\omega_2^-:=dx^1\wedge dx^3+dx^2\wedge dx^4$,
$\omega_3^-:=dx^1\wedge dx^4-dx^2\wedge dx^3$ be the basis for
anti-self dual $2$-forms on $\R^4$ introduced in $\S1.2$. We observe
that $\omega_1^-=d\beta_1^-$, $\omega_2^-=d\beta_2^-$ and
$\omega_3^-=d\beta_3^-$, with $\beta_1^-=x^1dx^2-x^3dx^4$,
$\beta_2^-=x^1dx^3+x^2dx^4$, $\beta_3^-=x^1dx^4-x^2dx^3$.

Let us define a family of connections on $B^4$ parameterized by real
$3\times 3$ matrices $\Asf:=(a_{ij})$ as follows
\begin{align}
\label{4.102}
\underline{B}_0(\Asf):=\underline{B}_{0,1}(\Asf)i+\underline{B}_{0,2}(\Asf)j+\underline{B}_{0,3}(\Asf)k\in
C^{\infty}(T^{\ast}B^4\otimes\so(3)),\notag\\
\text{ with\;
}\underline{B}_{0,l}(\Asf):=\underline{A}_{0,l}+a_{1l}\beta_1^-+a_{2l}\beta_2^-+a_{3l}\beta_3^-,\quad
1\le l\le 3,
\end{align}
and a family of  boundary connections
$B_0(\Asf)=\iota^{\ast}\underline{B}_0(\Asf)$. (Here, following the notation introduced in \cite{IM1},
$\;\underline{A}_{0}
:=\underline{A}_{0,1}i+\underline{A}_{0,2}j+\underline{A}_{0,3}k).$

Note that, for any choice of the matrix $\Asf$, the connection
$\underline{B}_0(\Asf)$ is a solution to $(\mathcal{D}_0)$ with
boundary value $B_0(\Asf)$. The corresponding matrix
$M(B_0(\Asf),p):= (m_{ij}(B_0(\Asf),p))_{1\le i,j\le 3}$ is computed explicitly as
follows:
\begin{align}
\label{4.103}
&m_{ij}(B_0(\Asf),p):=~\int_{B^4}\bigl((d\underline{B}_{0,j}(\Asf))^-,(dh_{p,i})^-\bigr)\notag\\
=&~m_{ij}(A_0,p)+a_{1j}\int_{B^4}\bigl(\omega_1^-,(dh_{p,i})^-\bigr)
+a_{2j}\int_{B^4}\bigl(\omega_2^-,(dh_{p,i})^-\bigr)
+a_{3j}\int_{B^4}\bigl(\omega_3^-,(dh_{p,i})^-\bigr).
\end{align}
Here, we recall that $h_p:=h_{p,1}i+h_{p,2}j+h_{p,3}k$, where the
components $h_{p,l}$ ($l=1,2,3$) are harmonic, and
$(dh_{p,l})^-:=(dh_{p,l})_1^-\omega_1^-+(dh_{p,l})_2^-\omega_2^-+(dh_{p,l})_3^-\omega_3^-$.
In this notation, $(\omega_k^-,(dh_{p,l})^-)=2(dh_{p,l})_k^-$. Since
the components $(dh_{p,l})_k^-$ are also harmonic on $B^4$, by the mean
value property for harmonic functions, one has
\begin{equation}
\label{4.104}
\int_{B^4}\bigl(\omega_k^-,(dh_{p,l})^-\bigr)
=2\int_{B^4}(dh_{p,l})_k^-=2|B^4|(dh_{p,l})_k^-(0)=\pi^2(dh_{p,l})_k^-(0),
\end{equation}
where $|B^4|=\frac{\pi^2}{2}$ is the Lebesgue measure of $B^4$.

\noindent Combining \eqref{4.103}, \eqref{4.104}, one obtains
\begin{equation}
\label{4.105}
 M(B_0(\Asf),p)=M(A_0,p)+\pi^2 H(p)\Asf,\quad\hbox{ with } H(p):=((dh_{p,l})_k^-(0))_{1\le l,k\le 3}\,.
\end{equation}

\noindent From \eqref{4.17} -- \eqref{4.19},  $H(p)$ can be written as
$H(p)=\left(\begin{array}{ccc}
\hsf_0(p)&\hsf_1(p)&\hsf_2(p)\\
-\hsf_1(p)&\hsf_0(p)&\hsf_3(p)\\
-\hsf_2(p)&-\hsf_3(p)&\hsf_0(p)
\end{array}\right)$,
with
\begin{align}
&\hsf_0(p)=\frac{1}{2}\Big(\frac{\partial\alpha_{p,1}}{\partial
x^1}+\frac{\partial\alpha_{p,2}}{\partial
x^2}+\frac{\partial\alpha_{p,3}}{\partial
x^3}+\frac{\partial\alpha_{p,4}}{\partial x^4}\Big)(0),\notag\\
&\hsf_1(p)=\frac{1}{2}\Big(-\frac{\partial\alpha_{p,1}}{\partial
x^4}+\frac{\partial\alpha_{p,2}}{\partial
x^3}-\frac{\partial\alpha_{p,3}}{\partial
x^2}+\frac{\partial\alpha_{p,4}}{\partial x^1}\Big)(0),\notag\\
&\hsf_2(p)=\frac{1}{2}\Big(\frac{\partial\alpha_{p,1}}{\partial
x^3}+\frac{\partial\alpha_{p,2}}{\partial
x^4}-\frac{\partial\alpha_{p,3}}{\partial
x^1}-\frac{\partial\alpha_{p,4}}{\partial x^2}\Big)(0),\notag\\
& \hsf_3(p)=\frac{1}{2}\Big(-\frac{\partial\alpha_{p,1}}{\partial
x^2}+\frac{\partial\alpha_{p,2}}{\partial
x^1}+\frac{\partial\alpha_{p,3}}{\partial
x^4}-\frac{\partial\alpha_{p,4}}{\partial x^3}\Big)(0).
\notag\end{align}
One has $\det H(p)=\hsf_0(p)(\hsf_0(p)^2+\hsf_1(p)^2+\hsf_2(p)^2+\hsf_3(p)^2)$ and
also, by \eqref{4.20}, $\hsf_0(p)>0$ for all $p\in B^4$. Thus, $H(p)$ is
non-singular for all $p\in B^4$.

\noindent This, in particular, implies the following result:

\newtheorem{proposition}{Proposition}[section]
\begin{proposition}
\label{P4.1} For any given $p_0\in B^4$, and any given $3\times 3$
matrix $M$, there exists always a $3\times 3$ matrix $\Asf$ such
that $M(B_0(\Asf),p_0)=M$. Thus, the different hypotheses of Theorem
3 can be obtained by choosing a suitable real matrix $\Asf$ (and
corresponding boundary data $B_0(\Asf)$).
\end{proposition}

Proposition \ref{P4.1}, combined with Theorem 3, yields many
examples of boundary values of the form $B_0(\Asf)$ such that
$(\mathcal{D}_{\epsilon})$ has multiple solutions in
$\mathcal{A}_{+1}(B_0(\Asf))$, as well as non-minimizing solutions.
Notice that, by elliptic regularity, the conditions in Theorems 1--3
are open conditions with respect to the boundary value (for example,
with respect to $H^{1/2}$-topology). The following proposition
holds.
\begin{proposition}
\label{P4.2}
For any given boundary connection $A_0\in
C^{\infty}(T^{\ast}\partial B^4\otimes\mathfrak{so}(3))$ and any
given $p_0\in B^4$, there exists an arbitrarily small perturbation
$\tilde{A}_0\in C^{\infty}(T^{\ast}\partial
B^4\otimes\mathfrak{so}(3))$ of $A_0$ such that the eigenvalues of
$M(\tilde{A}_0,p_0)^t M(\tilde{A}_0,p_0)$ satisfy
$\mu_1(\tilde{A}_0,p_0)>\mu_2(\tilde{A}_0,p_0)>\mu_3(\tilde{A}_0,p_0)>0$.
\end{proposition}
\textit{Proof:} We seek a boundary connection $\tilde{A}_0$ of the
form $B_0(\Asf)$ for some matrix $\Asf$ (as constructed in the
proof of Proposition \ref{P4.1}), with the desired requirements.

To this purpose, we first seek a $3\times 3$ matrix $X$ such that the
eigenvalues of $M(A_0,p_0)+X$ satisfy the requirements. We then
choose $\Asf$ such that $\pi^2H(p_0)\Asf=X$. For
$\tilde{A}_0=B_0(\Asf),$ the statement would then follow from
\eqref{4.105}.

Let us choose a matrix $Q\in SO(3)$ such that
$$Q^{-1}M(A_0,p_0)^tM(A_0,p_0)Q=\text{diag}(\mu_1(A_0,p_0),\mu_2(A_0,p_0),\mu_3(A_0,p_0)).$$
By taking an arbitrarily small perturbation, we may assume that
$M(A_0,p_0)$ is non-singular. For small $\mu>0$, we choose $X$ such
that $Q^{-1}M(A_0,p_0)^tXQ=\text{diag}(3\mu,2\mu,\mu)$ (take
$X=(M(A_0,p_0)^t)^{-1}Q\;\text{diag}\,(3\mu,2\mu,\mu)\;Q^{-1}$). We
then have
\begin{align}
\label{4.106}
&Q^{-1}(M(A_0,p_0)+X)^t(M(A_0,p_0)+X)Q\notag\\
=&~Q^{-1}M(A_0,p_0)^tM(A_0,p_0)Q+Q^{-1}M(A_0,p_0)^tXQ+Q^{-1}X^tM(A_0,p_0)Q+Q^{-1}X^tXQ\notag\\
=&~\text{diag}\,(\mu_1(A_0,p_0),\mu_2(A_0,p_0),\mu_3(A_0,p_0))+\text{diag}\,(6\mu,4\mu,2\mu)+O(\mu^2).
\end{align}
Thus, by taking $\Asf $ such that $\pi^2 H(p_0)\Asf=X$, we obtain
\begin{align}
&\mu_1(B_0(\Asf),p_0)=\mu_1(A_0,p_0)+6\mu+O(\mu^2),\notag\\
&\mu_2(B_0(\Asf),p_0)=\mu_2(A_0,p_0)+4\mu+O(\mu^2),\notag\\
&\mu_3(B_0(\Asf),p_0)=\mu_3(A_0,p_0)+2\mu+O(\mu^2).\notag
\end{align} These imply that there exists $\bar\mu >0$ such that, for
all $0<\mu<\bar\mu$, one has
$$\mu_1(B_0(\Asf),p_0)>\mu_2(B_0(\Asf),p_0)>\mu_3(B_0(\Asf),p_0)>0.$$
This completes the proof. \hfill$\Box$

\newtheorem{conjecture}{Conjecture}[section]
\begin{conjecture}
\label{R4.2} One could look for multiple solutions also in the
components $\mathcal{A}_k(A_0)$, for $k\ne \pm1$. We conjecture that
multiple solutions to $(\mathcal{D}_{\epsilon})$ exist in each
component $\mathcal{A}_k(A_0)$,  for small $\epsilon>0$, for a
rather general family of boundary values.
 For example, one may seek multiple solutions in
$\mathcal{A}_0(A_0)$ of the form $A=\underline
A_\epsilon\#\frac{1}{\epsilon}(\text{$1$-instanton})\#\frac{1}{\epsilon}(\text{$-1$-instanton})$.
This would require proofs similar to the ones established in this
paper, but with lengthier and more delicate calculations.
\end{conjecture}

\end{document}